\documentstyle[amscd,amssymb,verbatim,diagrams]{amsart}
\newcommand{\ad}{\operatorname{ad}}

\newcommand{\DD}{{\cal D}}

\newcommand{\splin}{\operatorname{sl}}

\newcommand{\gog}{{\frak g}}

\newcommand{\hra}{\hookrightarrow}
\newcommand{\lan}{\langle}
\newcommand{\ran}{\rangle}

\newcommand{\CC}{{\cal C}}

\newcommand{\Spec}{\operatorname{Spec}}

\renewcommand{\P}{{\Bbb P}}

\newcommand{\si}{\sigma}

\newcommand{\ga}{\gamma}
\newcommand{\Si}{\Sigma}
\newcommand{\de}{\delta}

\renewcommand{\ker}{\operatorname{ker}}

\newcommand{\bP}{{\bf P}}
\newcommand{\TCH}{{\cal T}{\operatorname{CH}}}
\newcommand{\TH}{{\cal T}H}
\numberwithin{equation}{section}

\newtheorem{thm}{Theorem}[section]
\newtheorem{prop}[thm]{Proposition}
\newtheorem{lem}[thm]{Lemma}
\newtheorem{cor}[thm]{Corollary}
\newenvironment{defi}{\vspace{3mm}\noindent
{\bf Definition.}}{\vspace{3mm}}
\newenvironment{rem}{\vspace{3mm}\noindent
{\bf Remark.}}{\vspace{3mm}}
\newenvironment{rems}{\vspace{3mm}
\noindent {\bf Remarks.}}{\vspace{3mm}}
\newenvironment{ex}{\vspace{3mm}\noindent
{\bf Example.}}{\vspace{3mm}}

\newcommand{\Pf}{\noindent {\it Proof}}
\newcommand{\id}{\operatorname{id}}

\newcommand{\pa}{\partial}

\renewcommand{\AA}{{\cal A}}

\newcommand{\PP}{{\cal P}}

\newcommand{\SS}{{\cal S}}
\newcommand{\LL}{{\cal L}}
\newcommand{\XX}{{\cal X}}
\newcommand{\YY}{{\cal Y}}
\newcommand{\RR}{{\cal R}}
\newcommand{\JJ}{{\cal J}}
\newcommand{\OO}{{\cal O}}

\renewcommand{\a}{\alpha}
\renewcommand{\b}{\beta}

\newcommand{\De}{\Delta}

\newcommand{\C}{{\Bbb C}}

\newcommand{\Z}{{\Bbb Z}}
\newcommand{\Q}{{\Bbb Q}}

\newcommand{\Ga}{\Gamma}

\newcommand{\wt}{\widetilde}
\newcommand{\ot}{\otimes}

\newcommand{\sub}{\subset}
\newcommand{\ed}{\qed\vspace{3mm}}

\newcommand{\HV}{{\cal HV}}
\newcommand{\CH}{\operatorname{CH}}
\newcommand{\TT}{{\cal T}}

\title{Algebraic cycles on the relative symmetric powers and on the
relative Jacobian of a family of curves. I}
\author{A. Polishchuk}
\address{Department of Mathematics, University of Oregon, Eugene, OR 97403}
\email{apolish@@uoregon.edu}
\thanks{This work was partially supported by the NSF grant DMS-0601034}

\begin{document}
\begin{abstract} In this paper we construct and study the actions of certain deformations of the Lie
algebra of Hamiltonians on the plane on the Chow groups (resp., cohomology) of the 
relative symmetric powers $\CC^{[\bullet]}$ and the relative Jacobian $\JJ$ of a family of curves 
$\CC/S$. As one of the applications, we show that 
in the case of a single curve $C$ this action induces a $\Z$-form of a Lefschetz $\splin_2$-action on the Chow groups of $C^{[N]}$. Another application gives a new grading on the ring $\CH_0(J)$
of $0$-cycles on the Jacobian $J$ of $C$ (with respect to the Pontryagin product) and equips it with an action of the Lie algebra of vector fields on the line. We also define the groups of tautological classes
in $\CH^*(\CC^{[\bullet]})$ and in $\CH^*(\JJ)$ and prove for them analogs of the properties established
in the case of the Jacobian of a single curve by Beauville in \cite{Bmain}. We show that the our
algebras of operators preserve the subrings of tautological cycles and act on them via some explicit differential operators. 
\end{abstract}
\maketitle

\bigskip

\centerline{\sc Introduction}

\bigskip

Let $\CC/S$ be a family of smooth projective curves over a smooth quasiprojective base $S$, 
and let $\CC^{[N]}$ denote the $N$th relative symmetric power of $\CC$ over $S$.
In this paper we construct and study the natural action of a certain modification of the Lie algebra of differential operators on the line on the direct sum of the Chow groups 
$$\CH^*(\CC^{[\bullet]}):=\bigoplus_{N\ge 0}\CH^*(\CC^{[N]}),$$
where $\CC^{[0]}=S$ (and on the similar direct sum of cohomology).
We also construct a related action of another algebra on $\CH^*(\JJ)$, where $\JJ/S$
is the corresponding relative Jacobian. These constructions are motivated by their potential use in the study of the Chow rings of Jacobians of curves, and in particular, in the study of the tautological subrings (see \cite{Bmain}, \cite{P-univ}, \cite{P-lie}). It was observed that the subalgebra generated by the standard cycles in the Jacobian of a curve depends in an interesting way on the corresponding point in the moduli space. Therefore, it is important to develop the corresponding calculus in the relative situation. On the other hand,
our main construction is reminiscent of the well known construction of the Heisenberg action on
cohomology of Hilbert schemes of surfaces (see \cite{Nak}, \cite{G}), and one can hope that there might
be a direct link between the two actions in the case of a family of curves lying on a surface.

Let us describe our main construction.
Let $\DD=\Z[t,\frac{d}{dt}]$ be the algebra of differential operators on the line.
Adjoin an independent variable $h$ and consider the subalgebra
$\DD_h\sub \DD\otimes\Z[h]$ generated over $\Z[h]$ by $t$ and by $h\frac{d}{dt}$.
We view it as a Lie algebra with the commutator
$$[D_1,D_2]_h=(D_1D_2-D_2D_1)/h.$$
Note that $\DD_h$ is a deformation of the Lie algebra $\wt{\HV}=\Z[t,p]$ of polynomial
Hamiltonians on the plane (equipped with the standard Poisson bracket).
Now for any supercommutative ring $A$ and an even element ${\bf a}_0\in A$ we define the Lie superalgebra 
$$\DD(A,{\bf a}_0):=\DD_h\otimes_{\Z[h]} A,$$
where the homomorphism $\Z[h]\to A$ sends $h$ to ${\bf a}_0$.
The (super)bracket  (resp., $\Z/2\Z$-grading) on $\DD(A,{\bf a}_0)$ is induced by the bracket on $\DD_h$ and
the product on $A$ (resp., by the $\Z/2\Z$-grading on $A$). More explicitly, 
$\DD(A,{\bf a}_0)$ is generated as an abelian group by the elements 
$$\bP_{m,k}(a):=t^m(h\frac{d}{dt})^k\ot a, \ \ a\in A.$$
The supercommutator is given by the formula
\begin{equation}\label{main-com-rel}
[\bP_{m,k}(a),\bP_{m',k'}(a')]= \sum_{i\ge 1}(-1)^{i-1}i!\cdot
\left({k\choose i}{m'\choose i}-{m\choose i}{k'\choose i}\right)
\bP_{m+m'-i,k+k'-i}(a\cdot a'\cdot {\bf a}_0^{i-1}).
\end{equation}
If $R\sub A$ is a subring then we can talk about an $R$-linear action of $\DD(A,{\bf a}_0)$ on an $R$-module (viewing $\DD(A,{\bf a}_0)$ as a Lie superalgebra over $R$). Our main construction gives a natural $\CH^*(S)$-linear (resp., $H^*(S)$-linear) action of $\DD(\CH^*(\CC),K)$ (resp.,
$\DD(H^*(\CC),cl(K))$) on $\CH^*(\CC^{[\bullet]})$ (resp., $\bigoplus H^*(\CC^{[N]})$),
where $K\in\CH^1(\CC)$ is the relative canonical class,
$cl(K)\in H^2(\CC)$ is the corresponding cohomology class (see Theorem \ref{action-thm} below).

For every integers $N\ge m\ge 0$ let us consider the morphism
$$s_{m,N}:\CC\times_S \CC^{[N-m]}\to \CC^{[N]},$$
sending a point $(p,D)\in \CC_s\times \CC_s^{[N-m]}$ to $mp+D$, where
$\CC_s\sub \CC$ is the fiber of our family over $s\in S$. Here
we identify points of $\CC_s^{[N]}$ with effective divisors of degree $N$ on $\CC_s$.  
Note that in the case $m=0$ this map is just the projection to $\CC^{[N]}$.
For a cycle $a\in\CH^*(\CC)$ and integers $m\ge 0$, $k\ge 0$, we consider the operator
$P_{m,k}(a)$ on $\CH^*(\CC^{[\bullet]})$ defined by the formula
\begin{equation}\label{main-oper-eq}
P_{m,k}(a)(x)=s_{m,N-k+m,*}(p_1^*a\cdot s_{k,N}^*x),
\end{equation}
where $x\in\CH^*(\CC^{[N]})$.
In the case $a=[\CC]$ we will simply write $P_{m,k}(\CC)$.
Note that in the case $N<k$ we have $P_{m,k}(a)(x)=0$.
If $a\in\CH^i(\CC)$ then
$P_{m,k}(a)$ sends $\CH^p(\CC^{[N]})$ to $\CH^{p+i+m-1}(\CC^{[N-k+m]})$.
%resp., $\CH_q(\CC^{[N]})$ to $\CH_{q-i-m+1}(\CC^{[N-k+m]})$.

\begin{thm}\label{action-thm}
(a) The map $\bP_{m,k}(a)\mapsto P_{m,k}(a)$ defines a $\CH^*(S)$-linear action of
$\DD(\CH^*(\CC),K)$ on $\CH^*(\CC^{[\bullet]})$,
where $K\in\CH^1(\CC)$ is the relative canonical class.
More precisely, the commutation relations \eqref{main-com-rel} for the operators $(P_{m,k}(a))$
(with ${\bf a}_0=K$) hold on the level of relative correspondences over 
$S$, i.e., they correspond to certain equalities in $\CH^*(\CC^{[\bullet]}\times_S \CC^{[\bullet]})$.

\noindent
(b) If we work over $\C$, the same construction defines an action of $\DD(H^*(\CC,\Z),cl(K))$ on 
$H^*(\CC^{[\bullet]})=\bigoplus_N H^*(\CC^{[N]},\Z)$.
\end{thm}

In the case of a trivial family $\CC=C\times S$ the above relations can be rewritten in a simpler form
(due to the fact that $\DD_h/h^2$ becomes a trivial deformation of $\wt{\HV}$).
Recall that the commutator in the Lie algebra $\wt{\HV}=\Z[x,p]$ of polynomial Hamiltonians on the plane is given by 
$$\{ x^mp^k, x^{m'}p^{k'} \}=(km'-mk') x^{m+m'-1}p^{k+k'-1}.$$

\begin{cor}\label{triv-base-cor}
Let $C$ be a smooth projective curve over a field $k$.
Choose a theta characteristic $\chi\in\CH^1(C)$ (so that $2\chi=K$) and set
$$L_{m,k}(a)=P_{m,k}(a)-mk P_{m-1,k-1}(p_1^*\chi\cdot a)$$ 
for $k\ge 0$, $m\ge 0$, $a\in\CH^*(C\times S)$. Then one has the following  relations:
$$[L_{m,k}(a),L_{m',k'}(a')]=(km'-mk')L_{m+m'-1,k+k'-1}(a\cdot a').$$
In other words, the map $x^mp^k\ot a\mapsto L_{m,k}(a)$ defines an action of
$\wt{\HV}\ot\CH^*(C\times S)$ on $\CH^*(C^{[\bullet]}\times S)$.
Similarly, we can define an action of $\wt{\HV}\ot H^*(C\times S)$ on $H^*(C^{[\bullet]}\times S)$.
We have $L_{0,0}([C\times S])=0$, so the operators $(L_{m,k}([C\times S]))$ define the action
of $\HV=\wt{\HV}/\Z$ on $\CH^*(C^{[\bullet]}\times S)$.
\end{cor}

\begin{rem} The operators $P_{n,0}(a)$ (resp., $P_{0,n}(a)$) for $n\ge 0$
are defined by the correspondences that are similar to those defining
Nakajima's operators $q_n(a)$ (resp., $q_{-n}(a)$) for the Hilbert schemes of points on a surface, 
where we use the notation of \cite{Lehn}. It is somewhat surprising that in the curve case
the Lie superalgebra generated by these operators is more complicated 
(the relations for $q_n(a)$ are simply those of the Heisenberg superalgebra). 
\end{rem}

Looking at the simplest of the above operators (such as $P_{m,1}(C)$ and $P_{0,1}([p])$, where 
$p\in C$ is a point) in the case $S=\Spec(k)$ we will derive the following result.

\begin{thm}\label{curve-thm} 
Let $C$ be a smooth projective curve of genus $g\ge 1$ over an algebraically closed field $k$ and let
$J$ be the Jacobian of $C$. Fix a point $p_0\in C(k)$ and consider the embedding $\iota:C\to J$
associated with $p_0$, so that $\iota(p_0)=0\in J(k)$. 
Let us denote by $I_C\sub\CH_0(J)$ the subgroup of classes represented by 
$0$-cycles of degree $0$ supported on $\iota(C)$. 

\noindent
(i) One has a direct sum decomposition 
$$\CH_0(J)=\Z\cdot [0]\oplus I_C\oplus I_C^{*2}\oplus\ldots\oplus I_C^{*g},$$
where $I_C^{*n}$ denotes the $n$th Pontryagin power of $I_C$. 
The associated filtration $(\bigoplus_{i\ge n}I_C^{*i})$ coincides with the 
standard filtration $(I^{*n})$, where $I\sub\CH_0(J)$ is the subgroup of cycles of degree zero.

\noindent
(ii) There exists a family of derivations $(\de_m)_{m\ge 1}$ of the graded algebra 
$$\CH_0(J)=\Z\oplus \bigoplus_{n=1}^g I_C^{*n}\simeq\bigoplus_{n=0}^g I^{*n}/I^{*(n+1)},$$
where the multiplication is given by the Pontryagin product, such that for every $x\in I_C$ one has
$$\de_m(x)=\sum_{i=0}^{m-1}(-1)^i{m\choose i}[m-i]_*x\in I_C^{*m}.$$
Equivalently, $\de_m|_{I_C}$ can be characterized by the property
$$\de_m([\iota(p)]-[0])=([\iota(p)]-[0])^{*m} \text{ for all }p\in C(k).$$
These derivations satisfy the commutation relations 
$$[\de_m,\de_{m'}]=(m'-m)\de_{m+m'-1},$$
i.e., they define an action of the Lie algebra of polynomial 
vector fields on the line vanishing at the origin by $t^m\frac{d}{dt}\mapsto\de_m$.
\end{thm}

We will show (see Remark 2 in the end of section \ref{first-sec}) 
that for a general curve of genus $\ge 3$ the decomposition in Theorem \ref{curve-thm}(i) 
is different from
the decomposition defined by Beauville (see \cite{B1}, p.254; \cite{B2}, Prop.~4). 

Theorems \ref{action-thm} and \ref{curve-thm} will be proved in section \ref{first-sec}. 
In sections \ref{Jac-sec} and \ref{div-sec}, that are somewhat more technical, 
we reprove (and generalize to the relative case) some known results using our algebra of operators. 

In section \ref{Jac-sec} we study the relation between our operators $(P_{m,k}(a))$ 
and the operators $(X_{n,k})_{n+k\ge 2}$ on $\CH(J)_{\Q}$,
where $J$ is the Jacobian of a curve $C$, constructed in \cite{P-lie}.
The latter family of operators satisfies the commutation relations of the 
Lie subalgebra $\HV'\sub\HV$ spanned over $\Z$ by the elements $x^np^k$ with 
$n+k\ge 2$ (it corresponds to Hamiltonian vector fields on the plane vanishing at the origin). 
It is natural to ask how this action is related to the one given by Corollary \ref{triv-base-cor} (say,
in terms of the push-forward map $\CH^*(C^{[\bullet]})_{\Q}\to\CH^*(J)_{\Q}$).
The relation turns out to be not quite straightforward. To work it out we
introduce another family of operators $(T_k(m,a))$, 
acting both on $\CH^*(C^{[\bullet]})$ and on $\CH^*(J)$, and
compatible with the push-forward map. In the case of $\CH^*(C^{[\bullet]})$ we find an explicit expression of $T_k(m,a)$ in terms of the operators $(P_{k,m}(a))$ (see Proposition \ref{T-P-prop}). 
On the other hand, in the
case of $\CH^*(J)_{\Q}$ we find that the operators $T_k(m,a)$ depend polynomially on $m$ and 
the corresponding coefficients are closely related to the operators $X_{n,k}$ considered in
\cite{P-lie}. Similar computations work in the case of the relative Jacobian
 $\JJ/S$ of a family of curves $\CC/S$. However, in the relative case the relations
 between the operators $X_{n,k}$ get deformed in an interesting way: the corresponding
(quadratic) algebra of operators acting on $\CH^*(\JJ)_{\Q}$ is not a universal enveloping algebra 
of a Lie algebra anymore (see \eqref{X-rel-eq}). We will study this algebra in detail elsewhere.

In section \ref{div-sec} we revisit algebraic Lefschetz $\splin_2$-actions for $C^{[N]}$ 
using our operators and study related questions of integrality.
The algebraic Lefschetz operators over $\Q$ are easily obtained from the action of
the algebra of (polynomial) differential operators in two variables $\DD_{t,u,\Q}$
on $\CH^*(\CC^{[\bullet]})_{\Q}$ generated by operators of the form $P_{10}(a)$ and $P_{01}(a')$
(see Corollary \ref{Heis-K-cor}).
After studying the divided powers of the operators $P_{n,0}(\CC)$ and $P_{0,n}(\CC)$ we construct an action of the divided powers subalgebra 
$\Z[t,u^{[\bullet]},\pa_t^{[\bullet]},\pa_u]\sub\DD_{t,u,\Q}$ on $\CH^*(\CC^{[\bullet]})$. Using this action
we generalize to the relative case 
the result of Collino \cite{Col2} on injectivity of the homomorphism 
$i_{N*}$ (resp., surjectivity of $i_N^*$), where $i_N:\CC^{[N-1]}\to\CC^{[N]}$ is the embedding
associated with $p_0$.
The corresponding $\splin_2$-action given by the operators $e=t\pa_u$, $f=u\pa_t$ and
$h=t\pa_t-u\pa_u$ preserves the Chow groups of the individual symmetric powers.
In the case when $S$ is a point we show that in this way we get a Lefschetz
$\splin_2$-triple for $C^{[N]}$ (see Theorem \ref{Lefschetz-thm}).
Working with divided powers allows us to reprove the fact 
(observed by del Ba\~{n}o in \cite{dB2}) 
that the hard Lefschetz isomorphism for $C^{[N]}$ holds over $\Z$ (see Corollary \ref{Lefschetz-cor}).

Finally, in section \ref{taut-sec} we define the groups of tautological classes in $\CH^*(\CC^{[\bullet]})$
and in $\CH^*(\JJ)$. For tautological classes in $\CH^*(\JJ)$ we establish the properties similar
to those obtained in the case $S=\Spec(k)$ by Beauville in \cite{Bmain}. We also show that tautological
subspaces are preserved by the operators constructed in this paper and by push-forward (resp.,
pull-back) associated with the relative Albanese maps $\CC^{[N]}\to\JJ$ (see Theorem \ref{taut-thm}). 
As an application of
our techniques we relate the modified diagonal classes in $\CH^{k-1}(C^{[k]})$ introduced by
Gross and Schoen in \cite{GS} to the pull-backs of some tautological classes on $J$ (see
Corollary \ref{pull-back-cor}). 

\vspace{3mm}

\noindent
{\it Notations and conventions}.

Throughout this paper we work with a family $\pi:\CC\to S$ of smooth projective curves of genus $g$,
where $S$ is smooth quasiprojective over a field $k$ (when we mention cohomology we assume
that $k=\C$). We denote by $\JJ/S$ the corresponding relative Jacobian.
In the case when $S$ is a point we denote $\CC$ (resp., $\JJ$) simply by $C$ (resp., $J$).
We will often use a natural product operation on $\CH^*(\CC^{[\bullet]})$ (resp., $\CH^*(\JJ)$)
called the {\it Pontryagin product}. It is defined using the map
$$\a_{m,n}:\CC^{[m]}\times_S \CC^{[n]}\to \CC^{[m+n]}:(D_1,D_2)\mapsto D_1+D_2.$$
by the formula
$$x*y=\a_{m,n,*}(p_1^*x\cdot p_2^*y)$$
for $x\in\CH^*(\CC^{[m]})$ and $y\in\CH^*(\CC^{[n]})$ (where $p_1$ and $p_2$ are the projections).
It is easy to see that this operation makes $\CH^*(\CC^{\bullet})$ into an associative commutative algebra over $\CH^*(S)=\CH^*(C^{[0]})$. The Pontryagin product on $\CH^*(\JJ)$ is defined similarly
using the addition map $\JJ\times_S\JJ\to\JJ$.
For every integer $m\in\Z$ we denote by $[m]:\JJ\to\JJ$ the corresponding map $\xi\mapsto m\xi$.
We usually fix a point $p_0\in\CC(S)$ and consider the corresponding embedding
$\iota:\CC\to\JJ$. We denote by $\LL$ the biextension on $\JJ\times_S\JJ$ corresponding to the autoduality of $\JJ$, normalized by the condition that $\LL|_{\CC\times_S\JJ}\simeq\PP_{\CC}$, where
$\PP_{\CC}$ is the Poincar\'e line bundle on $\CC\times_S\JJ$ trivialized over $p_0$ and over
the zero section of $\JJ$.

We often use results from Fulton's book \cite{Fulton}. 
We usually consider Chow groups
only for nonsingular varieties and use the upper grading (by codimension). 
For a cartesian diagram 
\begin{diagram}
X' &\rTo{} &Y'\\
\dTo{} & &\dTo{}\\
X &\rTo{f} & Y
\end{diagram}
where $f$ is a locally complete intersection morphism,
we denote by $f^!:\CH^*(Y')\to\CH^*(X')$ the refined Gysin map defined in section 6.6 of
\cite{Fulton}.

When we talk
about $0$-cycles we mean cycles of dimension zero and use the notation $\CH_0$.
Also, on one occasion 
in section \ref{Jac-sec} we also use Chow homology groups $\CH_*$ for a possibly singular scheme.
In the relative situation we view Chow groups of an $S$-scheme as a module over $\CH^*(S)$.
The analogs of our results for integral cohomology are based on the formalism developed in
\cite{FM} (that includes in particular Gysin maps $f_*:H^iX\to H^{i-2d}Y$ for proper locally complete
intersection morphisms $f:X\to Y$ such that $\dim X-\dim Y=d$).

The summation variables for which no range is given are supposed to be nonnegative integers.
The symbol $x^i$ for $i<0$ in algebraic formulas should be treated as zero. We denote divided powers
of a variable $x$ by $x^{[d]}=x^d/d!$.

\section{Cycles on symmetric powers}\label{first-sec}

In this section we will prove Theorems \ref{action-thm} and \ref{curve-thm}.
We start with computing some intersection products.
We will need to work with the closed embedding
$$t_{m,N}:\CC\times_S \CC^{[N-m]}\to \CC\times_S \CC^{[N]}:(p,D)\mapsto(p,mp+D).$$
Note that its composition with the projection to the second factor is the
map $s_{m,N}:\CC\times_S \CC^{[N-m]}\to \CC^{[N]}$ considered before.
We will denote by $\DD_N\sub \CC\times_S \CC^{[N]}$ the relative universal divisor (the image of
$t_{1,N}$).

\begin{lem}\label{mult-lem2}
Let us consider the cartesian diagram
\begin{diagram}
\Pi_{m,M,N} &\rTo{} & \CC^{[M]}\times_S \CC^{[N]}\\
\dTo{} & &\dTo{\a_{M,N}}\\
\CC\times_S \CC^{[M+N-m]} &\rTo{s_{m,M+N}}& \CC^{[M+N]}
\end{diagram}
For every decomposition $m=k+l$ we have a natural closed embedding
$$q_{k,l}:\CC\times_S \CC^{[M-k]}\times_S \CC^{[N-l]}\to \Pi_{m,M,N}:
(x,D_1,D_2)\mapsto (x,D_1+D_2,D_1+kx,D_2+lx),$$
where we view $\Pi_{m,M,N}$ as a subset of 
$\CC\times_S \CC^{[M+N-m]}\times_S\CC^{[M]}\times_S \CC^{[N]}$.
Then one has the following identity in $\CH^*(\Pi_{m,M,N})$:
$$\a_{M,N}^![\CC\times_S \CC^{[M+N-m]}]=
\sum_{k+l=m}{m\choose k} q_{k,l,*}[\CC\times_S \CC^{[M-k]}\times_S \CC^{[N-l]}].$$
\end{lem}

\Pf .
We have the following commutative diagram with cartesian squares
\begin{diagram}
\Pi_{m,M,N} &\rTo{} & \CC\times_S \CC^{[M]}\times_S \CC^{[N]} &\rTo{p_{23}} & 
\CC^{[M]}\times_S \CC^{[N]}\\
\dTo{} & & \dTo{\id\times\a_{M,N}} & &\dTo{\a_{M,N}}\\
\CC\times_S \CC^{[M+N-m]} &\rTo{t_{m,M+N}} & \CC\times_S \CC^{[M+N]} &\rTo{p_2}& \CC^{[M+N]}
\end{diagram}
Therefore, we have
$$\a_{M,N}^![\CC\times_S \CC^{[M+N-m]}]=(\id\times\a_{M,N})^![\CC\times_S \CC^{[M+N-m]}].$$
In the case $m=1$ the image of $t_{1,M+N}$ is exactly the universal divisor $\DD_{M+N}$.
By definition of the map $\a_{M,N}$, we have
\begin{equation}\label{div-pullback-eq}
(\id\times\a_{M,N})^*[\DD_{M+N}]=p_{12}^*[\DD_M]+p_{13}^*[\DD_N].
\end{equation}
This implies the required formula for $m=1$. The general case follows easily by induction in $m$.
\ed

\begin{lem}\label{diag-lem}
Consider the cartesian square
\begin{diagram}
\Si_{m,N} &\rTo{} & \CC\\
\dTo{} & &\dTo{\De_N}\\
\CC\times_S \CC^{[N-m]} &\rTo{s_{m,N}} &\CC^{[N]}
\end{diagram}
where $\De_N:\CC\to \CC^{[N]}$ is the relative diagonal embedding.
Note that there is a natural isomorphism $\Si_{m,N}\simeq \CC$ for $m>0$,
while $\Si_{0,N}\simeq \CC\times_S \CC$.
Then we have 
$$\De_N^!([\CC\times_S \CC^{[N]})=[\CC\times_S \CC]\in \CH^0(\CC\times_S \CC), \text{ and}$$
$$\De_N^!([\CC\times_S \CC^{[N-m]}])=(-1)^{m-1}m!{N\choose m}K^{m-1}\in\CH^{m-1}(\CC)$$
for $m\ge 1$.
\end{lem}

\Pf . The case $m=0$ is clear since $s_{0,N}$
is simply the projection $\CC\times_S \CC^{[N]}\to \CC^{[N]}$. In the case $m=1$
we should compute the intersection-product for the cartesian diagram
\begin{diagram}
\CC &\rTo{\De} & \CC\times_S \CC\\
\dTo{} & &\dTo{\id\times\De_N}\\
\CC\times_S \CC^{[N-1]} &\rTo{t_{1,N}} &\CC\times_S \CC^{[N]}
\end{diagram}
In other words, we have to compute the intersection of the universal divisor 
$\DD_N\sub \CC\times \CC^{[N]}$ with $\CC\times_S\De_N(\CC)$. We need to check that the corresponding multiplicity with the diagonal
$\De(\CC)\sub \CC\times_S \CC$ is equal to $N$. 
This is a local problem, so we can pick a local parameter $t$ along the fibers of $\CC\to S$
and think of $\DD_N$ as the set of pairs $(x,f(t))$, where $f$ is a unital polynomial of degree 
$N$ in $t$ such that $f(x)=0$. 
The diagonal embedding $\De_N$ sends a point $y$ to the polynomial $(t-y)^N$. Hence, the restriction of the equation $f(x)=0$ will have form $(x-y)^N$, which gives multiplicity $N$ with the diagonal $\De\sub \CC\times_S \CC$. 
The case of $m>1$ follows by induction: from the commutative diagram with cartesian squares
\begin{diagram}
\CC &\rTo{\id} & \CC &\rTo{\id} & \CC\\
\dTo{} & &\dTo{\id\times\De_{N-m+1}} & &\dTo{\De_N}\\
\CC\times_S \CC^{[N-m]} &\rTo{t_{1,N-m+1}} &\CC\times_S \CC^{[N-m+1]}&\rTo{s_{m-1,N}}&\CC^{[N]}
\end{diagram}
we see that
$$\De_N^!([\CC\times_S \CC^{[N-m]}])=(\id\times\De_{N-m+1})^!([\CC\times_S \CC^{[N-m]}]).$$
Hence, the step of induction follows from the previous computation (for $N-m+1$ instead of $N$)
along with the formula $\De^*([\De(\CC)])=K$.
\ed

\begin{lem}\label{inter-mult-lem} 
Consider the cartesian square
\begin{diagram}
Z_{m,k,N} &\rTo{} & \CC\times_S \CC^{[N-k]}\\
\dTo{} & & \dTo{s_{k,N}}\\
\CC\times_S \CC^{[N-m]} &\rTo{s_{m,N}} & \CC^{[N]}
\end{diagram}
We have natural closed embeddings
$$q^0:\CC\times_S \CC\times_S \CC^{[N-m-k]}\to Z_{m,k,N}: (x,x',D)\mapsto (x,D+kx',x',D+mx),$$
$$q^i:\CC\times_S \CC^{[N-m-k+i]}\to Z_{m,k,N}: (x,D)\mapsto (x,D+(k-i)x,x,D+(m-i)x),$$
where $1\le i\le\min(m,k)$ (we view $Z_{m,k,N}$ as a subset of 
$\CC\times_S \CC^{[N-m]}\times_S\CC\times_S \CC^{[N-k]}$).
Then we have the following formula for the intersection-product in the above diagram:
$$[\CC\times_S \CC^{[N-m]}]\cdot [\CC\times_S \CC^{[N-k]}]=
[q^0(\CC\times_S \CC\times_S \CC^{[N-m-k]})]+
\sum_{i\ge 1}(-1)^{i-1}i!{m\choose i}{k\choose i}q^i_*(K^{i-1}\times [C^{[N-m-k+i]}]).$$
\end{lem}

\Pf .
We can represent $s_{m,N}$ as the composition of $\De_m\times\id:\CC\times_S \CC^{[N-m]}\to 
\CC^{[m]}\times_S \CC^{[N-m]}$ followed by $\a_{m,N-m}:\CC^{[m]}\times_S \CC^{[N-m]}\to \CC^{[N]}$.
Therefore,
$$s_{m,N}^![\CC\times_S \CC^{[N-k]}]=(\De_m\times\id)^!\a_{m,N-m}^![\CC\times_S \CC^{[N-k]}].$$
Using Lemma \ref{mult-lem2} we obtain
$$s_{m,N}^![\CC\times_S \CC^{[N-k]}]=
\sum_{i+l=k}{k\choose i}z^{i,l}_*(\De_m\times\id)^![\CC\times_S \CC^{[m-i]}\times_S \CC^{[N-m-l]}],$$
where for $i+l=k$ we consider a closed subset $Z^{i,l}\stackrel{z^{i,l}}{\hra} Z=Z_{m,k,N}$
defined from the cartesian square
\begin{diagram}
Z^{i,l} &\rTo{} & \CC\times_S \CC^{[N-m]}\\
\dTo{} & & \dTo{\De_m\times\id}\\
\CC\times_S \CC^{[m-i]}\times_S \CC^{[N-m-l]} &\rTo{s^{i,l}} & \CC^{[m]}\times_S \CC^{[N-m]}
\end{diagram}
with
$$s^{i,l}(x,D_1,D_2)=(D_1+ix,D_2+lx).$$
Note that $s^{i,l}$ factors into the composition of the map
$$\CC\times_S \CC^{[m-i]}\times_S \CC^{[N-m-l]}\to \CC\times_S \CC^{[m-i]}\times_S \CC^{[N-m]}$$
induced by $t_{l,N-m}$ (identical on the second factor), followed by
$$s_{i,m}\times\id:\CC\times_S \CC^{[m-i]}\times_S \CC^{[N-m]}\to \CC^{[m]}\times_S \CC^{[N-m]}.$$
Now we can apply Lemma \ref{diag-lem}.
For $i=0$ we immediately get $Z^{0,k}\simeq \CC\times_S \CC\times_S \CC^{[N-m-k]}$ and
$$(\De_m\times\id)^![\CC\times_S \CC^{[m]}\times_S \CC^{[N-m-k]}]=
[\CC\times_S \CC\times_S \CC^{[N-m-k]}].$$
Similarly, for $i\ge 1$ we get
$Z^{i,k-i}=q^i(\CC\times_S \CC^{[N-m-k+i]})$ and
$$(\De_m\times\id)^![\CC\times_S \CC^{[m-i]}\times_S \CC^{[N-m-k+i]}]=(-1)^{i-1}i!{m\choose i}
p_1^*K^{i-1}\cdot [\CC\times_S \CC^{[N-m-k+i]}].$$
\ed

\noindent
{\it Proof of Theorem \ref{action-thm}.}
The operator $P_{k,m}(a)$ acting on $\CH^*(\CC^{[N]})$ is given by the relative correspondence
$f_{k,m*}(p_1^*(a))$, where 
$$f_{k,m}:\CC\times_S\CC^{[N-k]}\to\CC^{[N]}\times_S\CC^{[N-k+m]}:(p,D)\mapsto(kp+D,mp+D).$$
Therefore, to compute the correspondence inducing 
the composition $P_{k,m}(a)\circ P_{k',m'}(a')$ acting on $\CH^*(C^{[N]})$ we have to calculate
the push-forward to $\CC^{[N]}\times_S\CC^{[N'-k+m]}$ of 
the intersection product in the following diagram
\begin{diagram}
\CC\times_S \CC^{[N-k']}\times_S\CC^{[N'-k+m]} &&&&\CC^{[N]}\times_S\CC\times_S\CC^{[N'-k]}\\
&\rdTo{f_{k',m'}\times\id}&&\ldTo{\id\times f_{k,m}}\\
&&\CC^{[N]}\times_S\CC^{[N']}\times_S\CC^{[N'-k+m]}
\end{diagram}
multiplied with the pullbacks of $b$ and $a$, where $N'=N-k'+m'$,
It is easy to see that this intersection-product is exactly the one computed in 
Lemma \ref{inter-mult-lem} for $Z_{m',k,N'}\subset \CC\times_S\CC^{[N-k']}\times_S
\CC\times_S\CC^{[N'-k]}$. 
It follows that the above composition is given by the relative correspondence
$$(s_{k',N}\times s_{m,N'-k+m})_*(w\cdot p_1^*a'\cdot p_3^*a)\in\CH^*(\CC^{[N]}\times_S \CC^{[N'-k+m]}),$$
where $p_1, p_3: Z_{m',k,N'}\to \CC$ are the projections, and
$w\in\CH^*(Z_{m',k,N'})$ is the intersection-product computed in Lemma \ref{inter-mult-lem}.
When we substitute the formula for $w$ in the above equation and
subtract the similar expression for $P_{k',m'}(a')\circ P_{k,m}(a)$ we note that the first
terms (corresponding to the images of $q^0$) will cancel out due to the symmetry exchanging
the two factors $\CC$. The remaining terms will give the required formula for the commutator.
In the case of cohomology we have to work with the supercommutator since switching the order
of $a$ and $a'$ will introduce the standard sign.
\ed

\begin{cor}\label{Heis-K-cor} The operators $(P_{1,0}(a))$ and $(P_{0,1}(a))$ on 
$\CH^*(\CC^{[\bullet]})$ satisfy the following relations:
\begin{align*}
&[P_{1,0}(a),P_{1,0}(a')]=[P_{0,1}(a),P_{0,1}(a')]=0,\\
&[P_{0,1}(a),P_{1,0}(a')]=\lan a,a'\ran\cdot\id,
\end{align*}
where $\lan a,a'\ran=\pi_*(a\cdot a')\in\CH^*(S)$
(recall that we view $\CH^*(\CC^{[N]})$ as a $\CH^*(S)$-module
using the product with the pull-back under the projection $\CC^{[N]}\to S$).
In particular, if we are given a pair of divisor classes $\a,\b\in\CH^1(\CC)$ of nonzero relative degrees
$\deg(\a)$ and $\deg(\b)$ then there is an action of the algebra
$\DD_{t,u,\Q}=\Q[t,u,\pa_t,\pa_u]$ of differential operators in two variables on 
$\CH^*(\CC^{[\bullet]})_{\Q}$ such that
\begin{align*}
&t\mapsto P_{1,0}(\a)-\frac{\lan\a,\b\ran}{\deg(\b)}\cdot P_{1,0}(\CC) \ \
(\text{Pontryagin product with }\a-\frac{\lan\a,\b\ran}{\deg(\b)}\cdot[\CC])\\
&u\mapsto \frac{1}{\deg(\b)}P_{1,0}(\CC) \ \ (\text{Pontryagin product with } \frac{1}{\deg(\b)}[\CC]),\\ 
&\pa_t\mapsto \frac{1}{\deg(\a)}P_{0,1}(\CC),\\
&\pa_u\mapsto P_{0,1}(\b).
\end{align*}
For example, for $g\neq 1$ we can take $\a=\b=K$.
\end{cor}

\Pf . This follows from the relations of Theorem \ref{action-thm} together with the identity
$P_{0,0}(a)=\pi_*(a)\cdot\id$. 
\ed

\begin{rem} In the case $S=\Spec(\C)$
the cohomology $H^*(C^{[\bullet]},\Q)$ can be identified with the 
super-symmetric algebra of $H^*(C,\Q)$. Then the operators
$P_{1,0}(a)$ and $P_{0,1}(a)$ for $a\in H^*(C,\Q)$ are identified with the standard operators
on the super-symmetric algebra (products and contractions).
\end{rem}

\noindent
{\it Proof of Theorem \ref{curve-thm}.}
(i) Consider the Abel-Jacobi map $S:\CH_0(J)\to J(k): \sum m_i[a_i]\mapsto \sum m_ia_i$.
It is well known and easy to see that $S$ induces an isomorphism $I/I^{*2}\wt{\to} J(k)$ 
(see sec.0 of \cite{Bl}).
Since the composition $\CH_0(C)\stackrel{\iota_*}{\to}\CH_0(J)\stackrel{S}{\to} J(k)$ is the Abel-Jacobi
map for $C$ that induces an isomorphism of degree zero cycles with $J(k)$, we obtain
a decomposition
\begin{equation}\label{I-eq}
I=I_C\oplus I^{*2}.
\end{equation}
By taking the Pontryagin powers we immediately derive that
$$I^{*n}=I_C^{*n}+I^{*(n+1)}.$$
Since $I^{*(g+1)}=0$ by the result of Bloch (see \cite{Bl}), we deduce that
$$\CH_0(J)=\Z[0]+I_C+I_C^{*2}+\ldots+I_C^{*g}.$$
It remains to prove that this decomposition is direct, i.e., the summands are linearly independent.
The version of Roitman's theorem in 
arbitrary characterstic proved by Milne~\cite{Milne} implies that $I^{*2}\sub\ker(S)$ has no torsion. 
In view of \eqref{I-eq}, this shows that it is enough to prove our statement after tensoring with
$\Q$. A more direct way of getting our decomposition over $\Z$ will be outlined in Remark 3 after
Corollary \ref{module-cor}.

Let us consider the morphisms $C^{[N]}\to J$ induced by $\iota:C\to J$. 
It is easy to see that
the induced push-forward map $\si_*:\CH_0(C^{[\bullet]})\to\CH_0(J)$ is compatible with the
Pontryagin products. Also, since $\iota(p_0)=0\in J(k)$, we have $\si_*(x*[p_0])=\si_*(x)$ for
any $x\in\CH_0(C^{[\bullet]})$. Let $A_0(C)\sub\CH_0(C)$ denote the subgroup of classes of degree zero. Note that $\iota_*(A_0(C))=I_C$ and the push-forward
map $\si_*:\CH_0(C^{[g]})\to\CH_0(J)$ is an isomorphism (since $C^{[g]}\to J$ is birational).
Therefore, it suffices to establish the direct sum decomposition
\begin{equation}\label{AC-dec-eq}
\CH_0(C^{[g]})_{\Q}=
\Q\cdot [p_0]^{*g}\oplus A_0(C)_{\Q}*[p_0]^{*(g-1)}\oplus \ldots \oplus A_0(C)_{\Q}^{*(g-1)} *[p_0]\oplus 
A_0(C)_{\Q}^{*g},
\end{equation}
where the Pontryagin products are taken in $\CH_0(C^{[\bullet]})_{\Q}$.
To this end we will use the action of the algebra of differential operators $\DD_{t,\Q}$ in one variable
on $\CH^*(C^{[\bullet]})_{\Q}$ given by $t\mapsto P_{1,0}([p_0])$, 
$\frac{d}{dt}\mapsto P_{0,1}(C)$ (see Corollary \ref{Heis-K-cor}). Note that this
action preserves the subspace of $0$-cycles $\CH_0(C^{[\bullet]})_{\Q}$.
Since $\frac{d}{dt}$ acts locally nilpotently, we have a natural isomorphism of $\DD_{t,\Q}$-modules
\begin{equation}\label{D-mod-eq}
\CH_0(C^{[\bullet]})_{\Q}\simeq K_{\Q}[t],
\end{equation}
where $K=\ker(P_{0,1}(C))\cap\CH_0(C^{[\bullet]})$.
Furthermore, this isomorphism is compatible with gradings, where the grading on the right-hand side
is induced by the grading of $K$ and the rule $\deg(t)=1$. 
Next, we observe that for every $a\in\CH_0(C)$ we have the relation
$$[P_{0,1}(C),P_{1,0}(a)]=\deg(a)\cdot\id.$$
Hence, $P_{0,1}(C)$ commutes with the Pontryagin product with any $0$-cycle of degree zero on $C$.
Thus, we obtain 
$$A_0(C)^{*n}\sub K_n=\ker(P_{0,1}(C))\cap\CH_0(C^{[n]}) \text{ for }n\ge 1.$$
On the other hand, the algebra $\CH_0(C^{\bullet})_{\Q}$ is generated over $\Q=\CH_0(C^{[0]})_{\Q}$
by $\CH_0(C)_{\Q}=A_0(C)_{\Q}\oplus\Q\cdot [p_0]$.
Therefore, the natural map
$$\bigoplus_{n\ge 0}A_0(C)_{\Q}^{*n}[t]\to\CH_0(C^{[\bullet]})_{\Q}$$
is an isomorphism (where we set $A_0(C)_{\Q}^{*0}=\Q$). 
Looking at the grading component of degree $g$ we get the decomposition \eqref{AC-dec-eq}.

\noindent
(ii) Let us set $A=\Z\oplus A_0(C)\oplus A_0(C)^{*2}\oplus\ldots\oplus A_0(C)^{*g}$ (from (i)
we know that this algebra is isomorphic to $\CH_0(J)$) and consider
the natural homomorphism of algebras
\begin{equation}\label{A-hom-eq}
A[t]\to \CH_0(C^{[\bullet]})
\end{equation}
as in part (i). We claim that it is an isomorphism. Indeed, it is surjective, since $\CH_0(C^{[\bullet]})$
is generated by $[p_0]$ and by $A_0(C)$ as an algebra over $\Z$. Also, from part (i) we know that 
\eqref{A-hom-eq} is injective modulo torsion. It remains to show that it is injective
on the torsion subgroup. But the torsion in $A$ is contained in $A_0(C)$ (by Roitman's theorem), so the statement boils down to the fact that the natural map $A_0(C)t^N\to \CH_0(C^{[N]})$ is an embedding.
But this follows from the fact that its composition with the Abel-Jacobi map to $J(k)$ is an isomorphism.
Note that as in part (i) we could have avoided referring to Roitman's theorem and used
the divided powers instead (see Remark 3 after Corollary \ref{module-cor}).

Thus, we can view $(P_{m,1}(C))$ as operators on $A[t]$. For example, $P_{0,1}(C)$ acts by
$\frac{d}{dt}$.
For $m\ge 1$ and $a\in\CH_0(C)$ we have the relation 
$$[P_{m,1}(C),P_{1,0}(a)]=P_{m,0}(a).$$
This implies that for every point $p\in C(k)$ one has
$$P_{m,1}(C)([p]*x)=[p]*P_{m,1}(C)(x)+[p]^{*m}*x.$$
Since the classes $[p]$ generate our algebra, this implies that $P_{m,1}(C)$ is a derivation
of $A[t]$ characterized by 
$$P_{m,1}(C)([p])=[p]^{*m} \text{ for all }p\in C(k).$$
Setting $x_p=[p]-[p_0]\in A_0(C)$ we derive that $P_{m,1}(C)(t)=t^m$ and
$P_{m,1}(C)(x_p)=(x_p+t)^m-t^m$. Now let us define $\de_m$ as the following composition
$$A\to A[t]\stackrel{P_{m,1}(C)}{\rightarrow} A[t]\to A,$$
where the first map is the natural embedding and the last map is the evaluation at $t=0$.
Then $\de_m$ is a derivation of $A$ with the property $\de_m(x_p)=x_p^m$ for all $p\in C(k)$.
The commutation relations for $\de_m$ are easily checked on the generators $x_p$.
The formula for $\de_m|_{I_C}$ follows from the simple identity
$$([a]-[0])^{*m}=\sum_{i=0}^{m-1}(-1)^i{m\choose i}[m-i]_*([a]-[0])$$
in $\CH_0(J)$, where $m\ge 1$, $a\in J(k)$.
\ed

\begin{rems} 1. Note that $\de_1$ is just the grading derivation: it is equal to $n\id$ on the
grading component of degree $n$.
It is easy to see that under the identification \eqref{A-hom-eq} the operators
$P_{m,1}(C)$ are given by
$$P_{m,1}(C)=\de_m+mt\de_{m-1}+{m\choose 2}t^2\de_{m-2}+\ldots+mt^{m-1}\de_1+t^m\frac{d}{dt},$$
where $\de_m$ are extended to operators on $A[t]$ commuting with $t$.

\noindent 
2. It is natural to compare our decomposition of $\CH_0(J)=\CH^g(J)$ (tensored with $\Q$)
with the Beauville's decomposition
$$\CH^g(J)_{\Q}=\bigoplus_{s=0}^g\CH^g_s(J),$$
where $\CH^g_s(J)\sub\CH^g(J)_{\Q}$ is characterized by the condition $x\in\CH^g_s(J)$ if and only
if $[m]_*x=m^s x$ for all $m\in\Z$. The corresponding filtrations
$$\bigoplus_{s\ge n}\CH^g_s(J)=I^{*n}=\bigoplus_{s\ge n}I_C^*s$$
are the same (see \cite{B2}). However, the decompositions themselves are different.
Indeed, if they were the same we would have
$(I_C)_{\Q}\sub\CH^g_1(J)$ which would imply that $[2]_*([\iota(p)]-[0])=2[\iota(p)]-2[0]$ in
$\CH^g(J)_{\Q}$ for all $p\in C$. But this would mean that $([\iota(p)]-[0])^{*2}=0$ in $\CH^g(J)$, hence
$(p,p)-(p,p_0)-(p_0,p)+(p_0,p_0)$
is a torsion class in $\CH_0(C\times C)$ for all $p$, which is known not to be the case for a general
curve of genus $g\ge 3$ (see \cite{BV}, Prop. 3.2). 
\end{rems}

\section{Connection with cycles on the relative Jacobian}
\label{Jac-sec}

Assume that our family $\pi:\CC\to S$ is equipped with a section 
$p_0:S\to\CC$, and let $\si_N:\CC^{[N]}\to \JJ$
denote the corresponding map to the relative Jacobian of $\CC$
sending a divisor $D\in\CC_s^{[N]}$ to the class of the line bundle $\OO_{\CC_s}(D-Np_0)$.
We normalize the relative Poincar\'e line bundle $\PP_{\CC}$ on $\CC\times_S \JJ$ so 
that its pull-backs under $p_0\times\id_{\JJ}$ and under $\id_{\CC}\times e$ are trivial, where
$e:S\to\JJ$ is the zero section. Then we have the following equality in $\CH^1(\CC\times_S\CC^{[N]})$:
\begin{equation}\label{pull-back-P-b-eq}
(\id\times\si_N)^*c_1(\PP_{\CC})=[\DD_N]-Np_1^*[p_0]-p_2^*([\RR_N])-N\psi,
\end{equation}
where $\DD_N\sub \CC\times_S \CC^{[N]}$ is
the universal divisor, 
$$\RR_N=(p_0\times\id)^{-1}(\DD_N)=s_{1,N-1}(p_0\times \id)(\CC^{[N-1]})\sub \CC^{[N]}$$
is the divisor in $\CC^{[N]}$ associated with $p_0$, and
$$\psi=p_0^*K\in\CH^1(S)$$
(we view $\psi$ as a divisor class on any $S$-scheme via the pull-back).

We are going to introduce a new family of operators on $\CH^*(\JJ)$ and on
$\CH^*(\CC^{[\bullet]})$ that are compatible with respect to the push-forward map
$\si_*:\CH^*(\CC^{[\bullet]})\to\CH^*(\JJ)$ (that has $\si_{N*}$ as components).

It is convenient to consider a more general setup.
Let $\XX_{\bullet}=(\si_N:\XX_N\to\JJ)_{N\in\Z}$ be a family of proper $\JJ$-schemes equipped with
a collection of morphisms
$$s_N:\CC\times_S\XX_{N-1}\to \XX_N,$$
where $N\in\Z$, such that 

\noindent
(i) the diagram 
\begin{equation}
\begin{diagram}
\CC\times_S\XX_{N-1} & \rTo{s_N} & \XX_N\\
\dTo{\id\times\si_{N-1}} & &\dTo{\si_N}\\
\CC\times_S\JJ &\rTo{s_N}& \JJ
\end{diagram}
\end{equation}
is commutative, where the lower horizontal arrow is induced
by the map $\iota=\si_1:\CC\to\JJ$ and by the group law on the Jacobian;

\noindent
(ii) for each $N$ the map $\CC^m\times_S\XX_{N-m}\to\XX_N$ 
obtained from $(s_N)$ by iteration (where $\CC^m$ is the $m$th cartesian power of $\CC/S$), 
factors through a map
$\CC^{[m]}\times_S\XX_{N-m}\to\XX_N$.

Two main examples of the above situations are: $\XX_N=\CC^{[N]}$ for $N\ge 0$ 
(where $\si_N:\CC^{[N]}\to\JJ$ are associated with a point $p_0\in\CC(S)$, and $\XX_N=\emptyset$ for
$N<0$) and $\XX_N=\JJ$ for all $N\in\Z$. 
Another (singular) example is obtained by taking $\XX_N$
to be the image of the map $\CC^{[N]}\to\JJ$. 
%More generally, take $\XX_0$ to be any subvariety of $JJ$ and then...

Restricting the above maps $\CC^m\times_S\XX_{N-n}\to\XX_N$ to the diagonal in $\CC^m$
we get morphisms
$$s_{m,N}:\CC\times_S\XX_{N-m}\to\XX_N.$$ 
Now, let us define the operator $T_k(m,a)$ on 
$\CH_*(\XX_{\bullet})=\bigoplus_{N\in\Z}\CH_*(\XX_N)$, where $k\ge 0$, $m\ge 0$,
$a\in\CH^*(\CC)$, by the formula
$$T_k(m,a)(x)=s_{m,m+N,*}((\id\times\si_N)^*c_1(\PP_{\CC})^k\cdot p_1^*a\cdot p_2^*x),$$
where $x\in\CH_*(\XX_N)$, $p_1$ and $p_2$ are projections from the product 
$\CC\times_S \XX_N$ to its factors. In the case $a=[\CC]$ we will simply write $T_k(m,\CC)$.

Note that from the projection formula we get
\begin{equation}\label{T-k-0-eq}
T_k(0,a)(x)=\si_N^*\tau_k(a)\cdot x \ \text{ for }x\in\CH_*(\XX_N),
\end{equation}
where for $a\in\CH^*(\CC)$ and $k\ge 0$ we set 
\begin{equation}\label{tau-eq}
\tau_k(a)=p_{2*}(c_1(\PP_{\CC})^k\cdot p_1^*a)\in\CH^*(\JJ).
\end{equation}

Also, if $(f:\XX_N\to\YY_N)_{N\in\Z}$ is a morphism of two families as above then
it follows immediately from the definition that the above operators commute with the
push-forward map $f_*:\CH_*(\XX_{\bullet})\to\CH_*(\YY_{\bullet})$, i.e.,
$$T_k(m,a)\circ f_*=f_*\circ T_k(m,a).$$

\begin{thm}\label{relations-thm} 
One has the following relations between operators on $\CH_*(\XX_{\bullet})$: 
\begin{align*}
&\sum_{i\ge 0}\psi^i\cdot\left({k\choose i}m^{\prime i}
T_{k-i}(m,a)T_{k'}(m',a')-{k'\choose i}m^i T_{k'-i}(m',a')T_k(m,a)\right)=\\
&\sum_{i\ge 1}(-1)^{i-1}
\left({k\choose i}m^{\prime i}-{k'\choose i}m^i\right)T_{k+k'-i}(m+m',a\cdot a'\cdot
(K+2[p_0(S)])^{i-1})+\\
&\psi^{k'-1}m^{k'}p_0^*(a')T_k(m,a)-\psi^{k-1}m^{\prime k}p_0^*(a)T_{k'}(m',a')+\\
&\de_{k,0}\cdot\sum_{i\ge 1}{k'\choose i}m^i\psi^{i-1}p_0^*(a)T_{k'-i}(m',a')
-\de_{k',0}\cdot\sum_{i\ge 1}{k\choose i}m^{\prime i}\psi^{i-1}p_0^*(a') T_{k-i}(m,a)
\end{align*}
where $a,a'\in\CH^*(\CC)$, $k\ge 0$, $k'\ge 0$, $m\ge 0$, $m'\ge 0$. 
\end{thm}

Let us set $\ell_N=(\id\times\si_N)^*c_1(\PP_C)\in\CH^1(\CC\times_S \XX_N)$.
Also let us denote $\mu=\ell_1=(\id\times\si_1)^*c_1(\PP_C)\in\CH^1(\CC\times_S\CC)$.
Recall that $\PP_C$ is the pull-back of the biextension $\LL$ of $\JJ\times_S\JJ$
under the embedding $(\iota\times\id):\CC\times_S \JJ\to \JJ\times_S \JJ$
corresponding to point $p_0$. This implies the following isomorphism in
$CH^1(\CC\times_S\CC\times_S\XX_{N-m})$:
\begin{equation}\label{biext-eq}
(\id_{\CC}\times s_{m,N})^*\ell_N=m\cdot p_{12}^*\mu+p_{13}^*\ell_{N-m}.
\end{equation}

\begin{lem}\label{diag-lem2} 
One has the following identity in $\CH^*(\CC\times_S\CC)$ for $n\ge 1$:
$$\mu^n=(-\psi)^n+(-1)^n\psi^{n-1}\cdot
\left((p_0\times\id)_*[\CC]+(\id\times p_0)_*[\CC]\right)+(-1)^{n-1}
\sum_{i\ge 1}{n\choose i}\psi^{n-i}\cdot\De_*(K+2[p_0(S)])^{i-1}.$$
\end{lem}

\Pf . In the case $N=1$ the equality \eqref{pull-back-P-b-eq} gives
$$\mu=\De_*[\CC]-(p_0\times\id)_*[\CC]-(\id\times p_0)_*[\CC]-\psi\cdot [\CC\times_S\CC],$$
where $\De:\CC\to\CC\times_S\CC$ is the diagonal. Now the required identity is easily proved by induction in $n$ using the equalities 
$$\De^*\left(\De_*[\CC]-(p_0\times\id)_*[\CC]-(\id\times p_0)_*[\CC]\right)=-(K+2[p_0(S)]),$$
$$\mu\cdot (p_0\times\id)_*[\CC]=\mu\cdot (\id\times p_0)_*[\CC]=0.$$
\ed

\noindent
{\it Proof of Theorem \ref{relations-thm}.}
The composition $T_k(m,a)\circ T_{k'}(m',a')$ acting on $\CH_*(\XX_N)$ is given by
the operator
$$x\mapsto s_{m,m',*}\left((\id_{\CC}\times s_{m',m'+N})^*\ell_{m'+N}^{k}
p_{23}^*\ell_N^{k'}\cdot p_1^*(a)\cdot p_2^*(a')\cdot p_3^*(x)\right),$$
where $p_i$ ($i=1,2,3$) are the projections from the product $\CC\times_S\CC\times_S\XX_N$ to its factors, and $s_{m,m'}$ denotes the following composition
$$
\begin{diagram}
\CC\times_S\CC\times_S\XX_N &\rTo{\id_{\CC}\times s_{m',m'+N}}&\CC\times_S\XX_{m'+N}
&\rTo{s_{m,m+m'+N}}&\XX_{m+m'+N}.
\end{diagram}
$$ 
From \eqref{biext-eq} we get 
\begin{equation}\label{biext-pow-eq}
(\id_{\CC}\times s_{m',m'+N})^*\ell_{m'+N}^k=\sum_i {k\choose i}m^{\prime i} p_{12}^*\mu^i
\cdot p_{13}^*\ell_N^{k-i}.
\end{equation}
Let us set
$$S_{k,k';m,m'}(a,a')(x)=s_{m,m',*}\left(p_{12}^*\ell_N^k\cdot p_{23}^*\ell_N^{k'}\cdot p_1^*(a)\cdot p_2^*(a')\cdot p_3^*(x)\right).$$
Then using Lemma \ref{diag-lem2} and \eqref{biext-pow-eq} we derive
\begin{align*}
&T_k(m,a)\circ T_{k'}(m',a')=
\sum_{n\ge 0}(-m'\psi)^n{k\choose n}S_{k-n,k';m,m'}(a,a')+\\
&\sum_{n\ge 1,i\ge 1}(-1)^{n-1}m^{\prime n}\psi^{n-i}{k\choose n}{n\choose i}
T_{k+k'-n}(m+m',aa'(K+2[p_0(S)])^{i-1})+\\
&(-1)^km^{\prime k}\psi^{k-1}p_0^*(a)\cdot T_{k'}(m',a')+
\de_{k',0}\sum_{n\ge 1}(-1)^nm^{\prime n}\psi^{n-1}{k\choose n}p_0^*(a')\cdot T_{k-n}(m,a).
\end{align*}
From this we deduce that
\begin{align*}
&\sum_{p\ge 0}{k\choose p}(m'\psi)^p T_{k-p}(m,a)T_{k'}(m',a')=S_{k,k';m,m'}(a,a')+\\
&\sum_{i\ge 1}(-1)^{i-1}{k\choose i}m^{\prime i}T_{k+k'-i}(m+m',aa'(K+2[p_0(S)])^{i-1})
-m^{\prime k}\psi^{k-1}p_0^*(a)T_{k'}(m',a')\\
&-\de_{k',0}\cdot\sum_{i\ge 1}{k\choose i}m^{\prime i}\psi^{i-1}p_0^*(a')T_{k-i}(m,a).
\end{align*}
It remains to observe that condition (ii) imposed on $(\XX_N)$ implies that 
$$S_{k,k';m,m'}(a,a')=S_{k',k;m',m}(a',a).$$
Expressing both sides of this equality in terms of the operators $(T_k(m,a))$ we get the required
relation.
\ed

\begin{rems} 1. In the case of $\XX_N=\JJ$ the definition of the operators $T_k(m,a)$ 
can be extended to the case of arbitrary $m\in\Z$, so that the relations of Theorem \ref{relations-thm} 
still hold. Namely, we define the morphism $s_{m,N}$ for $m\in\Z$, so that it maps $(p,\xi)\in\CC_s\times\JJ_s$ to $m\iota(p)+\xi\in\JJ_s$.

\noindent
2. Similar relation holds if we replace Chow groups with cohomology 
provided we insert the standard sign $(-1)^{\deg(a)\deg(a')}$ whenever $a'$ goes before $a$.
\end{rems}

We are mostly interested in two cases: $(\XX_N=\CC^{[N]})$ and $(\XX_N=\JJ)$.
First, let us consider the case $(\XX_N=\CC^{[N]})$. In this case we will deduce the relation 
between the operators $(T_k(m,a))$ and $(P_{i,j}(a))$. We need one auxiliary result for this.
Let us denote by $Z_{m,N}\sub \CC\times_S \CC^{[N]}$ the image of the closed embedding $t_{m,N}$.
In other words, $Z_{m,N}$   
consists of $(p,D)\in\CC_s\times\CC_s^{[N]}$ such that $D-mp\ge 0$. Note that
$\DD_N:=Z_{1,N}$ is the universal divisor in $\CC\times_S \CC^{[N]}$.

\begin{lem}\label{Z-lem} 
One has the following equality in $\CH^m(\CC\times_S \CC^{[N]})$:
$$[Z_{m,N}]=[\DD_N]\cdot([\DD_N]+p_1^*K)\cdot\ldots\cdot([\DD_N]+(m-1)p_1^*K).$$
where $p_1:\CC\times_S \CC^{[N]}\to \CC$ is the projection.
Equivalently,
$$[\DD_N]^m=\sum_{i=0}^m (-1)^{m-i}S(m,i)\cdot p_1^*K^{m-i}\cdot [Z_{i,N}],$$
where $S(m,i)=\frac{1}{i!}\sum_{j=0}^i (-1)^j{i\choose j}(i-j)^m$ are the Stirling numbers of
the second kind.
\end{lem}

\Pf . In the case $m=1$ the equality is clear. The general case follows by induction in $m$ using the 
identity
\begin{equation}\label{D-res-eq}
[\DD_N]=t_{1,N+1}^*(\DD_{N+1}+p_1^*K)
\end{equation}
in $\CH^1(\CC\times_S \CC^{[N]})$. To obtain this identity one can start with the equality
\eqref{div-pullback-eq} (for $M=1$) in $\CH^1(\CC\times_S\CC\times_S\CC^{[N]})$ and then
apply the pull-back with respect to the diagonal embedding 
$\De\times\id:\CC\times_S\CC^{[N]}\to\CC\times_S\CC\times_S\CC^{[N]}$
(when the base is a point this was observed in Proposition 19.1 of \cite{P-av}).
\ed

\begin{prop}\label{T-P-prop} 
One has the following equality of operators on $\CH^*(\CC^{[\bullet]})$:
\begin{align*}
&T_k(m,a)=(-1)^kP_{m,0}(a\cdot [p_0])P_{1,1}([\CC])^k\psi^{k-1}+\\
&\sum_{i+n+j=k}(-1)^{n+j}{k\choose j}S(i+n,i)
P_{i+m,i}(a\cdot K^n)P_{1,1}([p_0]+\psi)^j.
\end{align*}
\end{prop}

\Pf . Since $[Z_{k,N}]=t_{k,N,*}([\CC\times_S \CC^{[N-k]}])$,
we can rewrite $P_{k+m,k}(a)$ in the form similar to
that of $T_{k}(m,a)$:
\begin{equation}\label{P-k-m-Z-eq}
P_{k+m,k}(a)(x)=s_{m,N,*}([Z_{k,N}]\cdot p_1^*a\cdot p_2^*x),
\end{equation}
where $x\in\CH^*(\CC^{[N]})$.
Let us use the following shorthand notation for divisors on $\CC\times_S\CC^{[N]}$:
$\DD=\DD_N$, $\RR=p_2^*\RR_N$, $K=p_1^*K$,
$[p_0]=p_1^*[p_0]$. Then we can write \eqref{pull-back-P-b-eq} 
as
$$(\id\times\si_N)^*c_1(\PP_{\CC})=\DD-\RR-N[p_0]-N\psi.$$
Note that we have the following relations in $\CH^2(\CC\times_S\CC^{[N]})$:
$$(\DD-\RR)\cdot[p_0]=0, \ \ [p_0]^2=-\psi\cdot[p_0].$$
It follows that for $j\ge 1$ one has $(\DD-\RR-N\psi)^i\cdot [p_0]^j=N^i(-\psi)^{i+j-1}\cdot [p_0]$,
so we derive
$$(\id\times\si_N)^*c_1(\PP_{\CC})^k=(\DD-\RR-N([p_0]+\psi))^k=(\DD-\RR-N\psi)^k+(-N)^k[p_0]\psi^{k-1}.
$$
Therefore, using Lemma \ref{Z-lem} we find
$$(\id\times\si_N)^*c_1(\PP_{\CC})^k=(-N)^k[p_0]\psi^{k-1}+
\sum_{i+n+j=k}(-1)^{n+j}{k\choose j}S(i+n,i)[Z_{i,N}]\cdot K^n\cdot (\RR+N\psi)^j.$$
Taking into the account \eqref{P-k-m-Z-eq}, we get
\begin{align*}
&T_k(m,a)(x)=(-N)^kP_{m,0}(a\cdot[p_0])(x)\cdot\psi^{k-1}+\\
&\sum_{i+n+j=k}(-1)^{n+j}{k\choose j}S(i+n,i)
P_{i+m,i}(a\cdot K^n)((\RR+N\psi)^j\cdot x),
\end{align*}
where $x\in\CH^*(\CC^{[N]})$.
It remains to use the equalities
$$P_{1,1}([\CC])(x)=Nx,$$
$$P_{1,1}([p_0])(x)=\RR\cdot x$$
for $x\in\CH^*(\CC^{[N]})$.
\ed

Now let us specialize to the case $\XX_N=\JJ$. 
In this case we can relate the operators $T_k(m,a)$
to the operators considered in \cite{P-lie} (in the case $S=\Spec(k)$).
Following \cite{P-lie} let us define the operators $A_k(\a)$ on $\CH^*(\JJ)$ for $\a\in\CH^*(\JJ)$
and $k\ge 0$ by
$$A_k(\a)(x)=(p_1+p_2)_*(c_1(\LL)^k\cdot p_1^*\a\cdot p_2^*x),$$
where $p_1$ and $p_2$ are projections from the product $\JJ\times_S \JJ$ to its factors.

\begin{lem}\label{A-T-lem1} 
For $a\in\CH^*(\CC)$ one has 
$A_k([m]_*\iota_*a)=m^k T_k(m,a)$ (with the convention that $0^0=1$),
where $[m]:\JJ\to \JJ:\xi\to m\xi$.
\end{lem}

\Pf . We have
\begin{align*}
&A_k([m]_*\iota_*a)=(mp_1+p_2)_*\left(([m]\times\id)^*c_1(\LL)^k\cdot p_1^*(\iota_*a)\cdot p_2^*x\right)=\\
&m^k(mp_1+p_2)_*\left(c_1(\LL)^k\cdot (\iota\times\id_{\JJ})_*(p_1^*a)\cdot p_2^*x\right).
\end{align*}
Now the result follows immediately from the isomorphism $\LL|_{\CC\times_S \JJ}\simeq\PP_{\CC}$.
\ed

Working with rational coefficients we can
consider the decomposition
$$\CH^*(\JJ)_{\Q}=\bigoplus_{i=0}^{2g}\CH^*(\JJ)_i,$$
where $[m]_*x=m^ix$ for $x\in\CH^*(\JJ)_i$ (see \cite{DM} Thm. 3.1).
It follows that for fixed $k$ and $\a\in\CH^*(\JJ)$ the operator valued function
$m\to A_k([m]_*\a)$ is a polynomial of degree $\le 2g$. By Lemma \ref{A-T-lem1} the same is true for
$T_k(m,a)$, where $a\in\CH^*(\CC)$ (more precisely, it is a polynomial in $m$ of degree $\le 2g-k$). Therefore, we can write
$$T_k(m,a)=\sum_{n=0}^{2g-k} \frac{m^n}{n!} \wt{X}_{n,k}(a)$$
for some operators $\wt{X}_{n,k}(a)$ on $\CH^*(\JJ)_{\Q}$.
Then the relations of Theorem \ref{relations-thm}
are equivalent to the following relations for $(\wt{X}_{n,k}(a))$:
\begin{equation}\label{X-rel-eq}
\begin{array}{l}
\sum_{i\ge 0}\psi^i\cdot i!\left({k\choose i}{n'\choose i} \wt{X}_{n,k-i}(a)\wt{X}_{n'-i,k'}(a')-
{k'\choose i}{n\choose i} \wt{X}_{n',k'-i}(a')\wt{X}_{n-i,k}(a)\right)=\\
\sum_{i\ge 1}(-1)^{i-1}i!\left({k\choose i}{n'\choose i}-
{k'\choose i}{n\choose i}\right) \wt{X}_{n+n'-i,k+k'-i}(aa'(K+2[p_0(S)])^{i-1})+\\
\de_{n',0}p_0^*(a')\psi^{k'-1}k'!{n\choose k'}\wt{X}_{n-k',k}(a)-
\de_{n,0}p_0^*(a)\psi^{k-1}k!{n'\choose k}\wt{X}_{n'-k,k'}(a')+\\
\de_{k,0}p_0^*(a)\psi^{n-1}n!{k'\choose n}\wt{X}_{n',k'-n}(a')-
\de_{k',0}p_0^*(a')\psi^{n'-1}n'!{k\choose n'}\wt{X}_{n,k-n'}(a).
\end{array}
\end{equation}
We will denote $\wt{X}_{n,k}([\CC])$ simply by $\wt{X}_{n,k}(\CC)$.
In the case $S=\Spec(k)$ the above relations are essentially equivalent to those of
Theorem 2.6 of \cite{P-lie}. Recall that in \cite{P-lie} we showed that
the operators $\wt{X}_{n,k}(C)-nk\wt{X}_{n-1,k-1}(K/2+[p_0])$ 
satisfy the commutation relations of the Lie 
algebra $\HV'$ and calculated their Fourier transform. We are going to 
present a similar computation in the relative case (i.e., when $S$ is arbitrary).
We refer to \cite{DM} for the basic properties of the Fourier transform on cycles over
abelian schemes (originally introduced in \cite{M} and studied in \cite{B1} and \cite{B2}).

\begin{thm}\label{four-thm} 
(i) The operators
\begin{equation}
e=\frac{1}{2}\wt{X}_{0,2}(\CC),\ \ 
f=-\frac{1}{2}\wt{X}_{2,0}(\CC),\ \ 
h=-\wt{X}_{1,1}(\CC)+g\cdot\id
\end{equation}
on $\CH(\JJ)_{\Q}$ define an action of the Lie algebra $\splin_2$.

\noindent (ii) Let us set for $a\in\CH^*(\CC)$, $n\ge 0$, $k\ge 0$,
$$X_{n,k}(a)=\sum_{i\ge 0}(-1)^i i!{n\choose i}{k\choose i}\wt{X}_{n-i,k-i}(a\eta^i),$$
where $\eta:=K/2+[p_0(S)]+\psi/2\in\CH^1(\CC)_{\Q}$.
Consider the Fourier transform defined by
$$F:\CH^*(\JJ)_{\Q}\to\CH^*(\JJ)_{\Q}:x\mapsto p_{2*}(\exp(c_1(\LL))\cdot p_1^*x).$$
Then one has
$$FX_{n,k}(a)F^{-1}=(-1)^kX_{k,n}(a).$$
$$[e,X_{n,k}(a)]=nX_{n-1,k+1}, \ \ [f,X_{n,k}(a)]=kX_{n+1,k-1}(a), \ \ [h,X_{n,k}(a)]=(k-n)X_{n,k}(a).$$
\end{thm}

\begin{rem} The $\splin_2$-action of Theorem \ref{four-thm} differs only by a sign from 
the relative Lefschetz action associated with the relatively ample class $-\tau_2(\CC)/2$ on $\JJ$
(see \cite{K}).
\end{rem}

We start with some preliminary statements.

\begin{lem}\label{four-lem} 
One has $F\wt{X}_{n,0}(a)F^{-1}=\wt{X}_{0,n}(a)$ for every $a\in\CH^*(\CC)$, $n\ge 0$.
\end{lem}

\Pf . By definition, $\wt{X}_{n,0}(a)/n!$ is the Pontryagin product with $a_n\in\CH^*(\JJ)_{\Q}$
where $[m]_*\iota_*a=\sum_{i\ge 0}m^ia_i$ for all $m\in\Z$. On the other hand,
$\wt{X}_{0,n}(a)=T_n(0,a)$ 
is the usual product with $\tau_n(a)$ (see \eqref{T-k-0-eq} and \eqref{tau-eq}).
Since $F(x*y)=F(x)\cdot F(y)$, it remains to show that 
\begin{equation}\label{F-a-eq}
F(a_n)=\frac{1}{n!}\tau_n(a)=p_{2*}(\frac{c_1(\LL)^n}{n!}\cdot p_1^*\iota_*a),
\end{equation}
where $p_1$ and $p_2$ are the projections of the product $\JJ\times_S\JJ$ on its factors.
This fact is well known but we will give the proof since it is very short. We have
$$\sum_{i\ge 0}m^i F(a_i)=F([m]_*\iota_*a)=[m]^*F(\iota_*a)=[m]^*p_{2*}(\exp(c_1(\LL))\cdot
p_1^*\iota_*a)=p_{2*}((\id\times[m])^*\exp(c_1(\LL))\cdot p_1^*\iota_*a).$$
Using the identity $(\id\times [m])^*c_1(\LL)=m c_1(\LL)$, we see that this is equal
to $\sum_{i\ge 0}m^i\tau_i(a)/i!$. Now \eqref{F-a-eq} is obtained by equating the coefficients with $m^n$. 
\ed

\begin{lem}\label{sl2-lem}
(i) One has $\wt{X}_{0,i}(\CC)=\wt{X}_{i,0}(\CC)=0$ for $i\le 1$.

\noindent (ii) One has
$$[e,\wt{X}_{n,k}(a)]=n\wt{X}_{n-1,k+1}(a)-n(n-1)\wt{X}_{n-2,k}(a\cdot\eta),$$
$$[f,\wt{X}_{n,k}(a)]=k\wt{X}_{n+1,k-1}(a)-k(k-1)\wt{X}_{n,k-2}(a\cdot\eta).$$
\end{lem}

\Pf . (i) The operator $\wt{X}_{0,0}(\CC)=T_0(0,\CC)$ is the product with the pull-back of
$\pi_*[\CC]=0$, where $\pi:\CC\to S$ is the projection. On the other hand, the operator
$\wt{X}_{0,1}(\CC)=T_1(0,\CC)$ is the product with $\tau_1(\CC)=p_{2*}(c_1(\PP_{\CC}))$,
where $p_2:\CC\times_S\JJ\to\JJ$ is the projection. Since $c_1(\PP_{\CC})$ is the divisor
in $\CC\times_S\JJ$ that has degree zero on every fiber of $p_2$, we get $\tau_1(\CC)=0$.
Now the vanishing of $\wt{X}_{1,0}(\CC)$ follows from Lemma \ref{four-lem}.

\noindent
(ii) This is an immediate consequence of relations \eqref{X-rel-eq} and of part (i).
\ed

\noindent
{\it Proof of Theorem \ref{four-thm}}.
(i) These relations follow from Lemma \ref{sl2-lem} and from the observation that
$\wt{X}_{0,0}(\eta)=g\cdot\id$ (since it is the product with the pull-back of $\pi_*(\eta)=g\cdot [S]$).

\noindent
(ii) From Lemma \ref{four-lem} we get 
\begin{equation}\label{F-f-e-eq}
F f F^{-1}=-e. 
\end{equation}
Recall that $-f$ is the Pontryagin product with
the class $a_2$ on $\JJ$ defined by $[m]_*[\CC]=\sum_i m^i a_i$. Replacing $m$ by $-m$ we see
that $[-1]_*a_2=a_2$. Hence, $f$ commutes with $[-1]^*$. Now the identity 
$F^2=(-1)^g[-1]_*$ (see \cite{DM}, Cor.2.22) implies that $F^2$ commutes with $f$. Therefore, from
\eqref{F-f-e-eq} we get
$$FeF^{-1}=-f.$$
On the other hand, from
Lemma \ref{sl2-lem}(ii) we deduce by induction that
\begin{equation}\label{ad-eq}
\begin{array}{l}
\ad(f)^k \frac{\wt{X}_{0,n+k}(a)}{(n+k)!}=\frac{X_{k,n}(a)}{n!},\\
\ad(e)^k \frac{\wt{X}_{n+k,0}(a)}{(n+k)!}=\frac{X_{n,k}(a)}{n!}.
\end{array}
\end{equation}
Now, combining all of this with Lemma \ref{four-lem} we get
$$F\frac{X_{n,k}(a)}{n!}F^{-1}=\ad(-f)^k(F\frac{\wt{X}_{n+k,0}(a)}{(n+k)!}F^{-1})=
(-1)^k\ad(f)^k\frac{\wt{X}_{0,n+k}(a)}{(n+k)!}=(-1)^k\frac{X_{k,n}(a)}{n!}$$
as required. The formulas for $[e,X_{n,k}(a)]$ and $[f,X_{n,k}(a)]$ follow immediately
from \eqref{ad-eq}. 
\ed

\section{Divided powers}
\label{div-sec}

In this section we construct and study divided powers of operators 
$P_{n,0}(\CC)$ and $P_{0,n}(\CC)$. Then we use them to define an action of a certain
$\Z$-form of an algebra of differential operators in two variables on $\CH^*(\CC^{[\bullet]})$.
We will also construct a $\Z$-version of the Lefschetz $\splin_2$-action on $\CH^*(C^{[N]})$
(see Theorem \ref{Lefschetz-thm}).

Let us start with divided powers of $P_{n,0}(\CC)$. By definition, we have
$$P_{n,0}(\CC)(x)=\de_n * x,$$
where $\de_n=\De_{n*}([\CC])\sub\CH_1(\CC^{[n]})$ is the class of the diagonal.
Thus, we can set
\begin{equation}\label{div-eq1}
P_{n,0}(\CC)^{[d]}(x)=\de_n^{[d]}*x,
\end{equation}
where 
$$\de_n^{[d]}=\De_{n*}^{[d]}([\CC^{[d]}])\in\CH_d(\CC^{[nd]}),$$
$$\De_n^{[d]}:\CC^{[d]}\to \CC^{[nd]}: D\mapsto nD.$$ 
Note that $d!\de_n^{[d]}=\de_n^{*d}$ --- the $d$th power of $\de_n$ with the respect to the
Pontryagin product. Hence, $d!P_{n,0}(\CC)^{[d]}=P_{n,0}(\CC)^d$.

To describe the divided powers of $P_{0,n}(\CC)$ let us introduce a new binary operation on 
$\CH^*(\CC^{[\bullet]})$ as follows. For $a\in\CH^*(\CC^{[k]})$ and $x\in\CH^*(\CC^{[N]})$ set
$$i_a(x)=p_{2*}(p_1^*a\cdot\a_{k,N-k}^*x)\in\CH^*(\CC^{[N-k]}),$$
where $p_1,p_2$ are projections of the product $\CC^{[k]}\times_S \CC^{[N-k]}$ to its factors.
Then it is easy to see that 
$$P_{0,n}(\CC)(x)=i_{\de_n}(x).$$
Also, it is straightforward to check that
$$i_{a*b}=i_a\circ i_b.$$
Thus, it is natural to set
\begin{equation}\label{div-eq2}
P_{0,n}(\CC)^{[d]}(x)=i_{\de_n^{[d]}}(x),
\end{equation}
so that we have $d!P_{0,n}(\CC)^{[d]}=P_{0,n}(\CC)^d$.

Below we use the notation from the Introduction.
Let $A$ be a supercommutative algebra $A$ with a unit and a distinguished
even element ${\bf a}_0\in A$. We are going to define two extensions of the universal enveloping algebra
of $\DD(A,{\bf a}_0)$ by adding two families of divided powers.

\begin{lem}\label{div-sum-lem}
(i) Let $\gog$ be a Lie algebra over $\Z$. Then for $x,y\in\gog$, the following relations hold in
$U(\gog)$:
$$x^{d}y=\sum_{i=0}^d {d\choose i}(\ad x)^{i}(y)x^{d-i},$$
$$yx^{d}=\sum_{i=0}^d {d\choose i}x^{d-i}(-\ad x)^{i}(y) $$
for all $d\ge 1$.

\noindent (ii) The following relations hold in $U(\DD(A,{\bf a}_0))$:
$$\bP_{m,k}(a)\bP_{n,0}(1)^{d}=\sum_{j\le i}(-1)^{i-j}\frac{i!d!}{j!(d-j)!}{k\choose i}A_j(i,n)\bP_{n,0}(1)^{d-j}
\bP_{m+nj-i,k-i}(a{\bf a}_0^{i-j}),$$
$$\bP_{0,n}(1)^{d}\bP_{m,k}(a)=\sum_{j\le i}(-1)^{i-j}\frac{i!d!}{j!(d-j)!}{m\choose i}A_j(i,n)
\bP_{m-i,k+nj-i}(a{\bf a}_0^{i-j})\bP_{0,n}(1)^{d-j},$$
where
$$A_d(i,n)=\sum_{i_1+\ldots+i_d=i,i_s\ge1}{n\choose i_1}\ldots {n\choose i_d}.$$
We also use the convention $x^i=0$ for $i<0$ (so that $j\le d$ in both sums).
\end{lem}

\Pf . (i) This is easily checked by induction in $d$.

\noindent (ii) Using the commutation relations in $\DD(A,{\bf a}_0)$ one can check by induction in $d$ that
$$\ad(\bP_{0,n}(1))^d(\bP_{m,k}(a))=\sum_{i\ge d}(-1)^{i-d}i!{m\choose i}A_d(i,n)
\bP_{m-i,k+nd-i}(a\cdot {\bf a}_0^{i-d}),$$
$$\ad(-\bP_{n,0}(1))^d(\bP_{m,k}(a))=\sum_{i\ge d}(-1)^{i-d}i!{k\choose i}A_d(i,n)
\bP_{m+nd-i,k-i}(a\cdot {\bf a}_0^{i-d}).$$
Now the required relations follow from (i).
\ed

\begin{defi}
Let us denote by $\wt{U}_1(A,{\bf a}_0)$ (resp., $\wt{U}_2(A,{\bf a}_0)$)
the superalgebra over $\Z$ with generators
$(\bP_{m,k}(a))$, $m\ge 0, k\ge 0$, depending additively on $a\in A$, and $(\bP_{n,0}(1)^{[d]})$
(resp., $(\bP_{0,n}(1)^{[d]})$), $n\ge 1$,
$d\ge 0$, subject to the following relations:

\noindent
(i) the supercommutator relations of $\DD(A,{\bf a}_0)$ for $(\bP_{m,k}(a))$;

\noindent
(ii) $d!\bP_{n,0}(1)^{[d]}=\bP_{n,0}(1)^d$ (resp., $d!\bP_{0,n}(1)^{[d]}=\bP_{0,n}(1)^d$); \\
$\bP_{n,0}(1)^{[d_1]}\bP_{n,0}(1)^{[d_2]}={d_1+d_2\choose d_1}\bP_{n,0}(1)^{[d_1+d_2]}$
(resp.,
$\bP_{0,n}(1)^{[d_1]}\bP_{0,n}(1)^{[d_2]}={d_1+d_2\choose d_1}\bP_{0,n}(1)^{[d_1+d_2]}$);

\noindent
(iii)
$\bP_{m,k}(a)\bP_{n,0}(1)^{[d]}=\sum_{j\le i}(-1)^{i-j}\frac{i!}{j!}{k\choose i}A_j(i,n)\bP_{n,0}(1)^{[d-j]}
\bP_{m+nj-i,k-i}(a{\bf a}_0^{i-j})$, 
$$(\text{resp., } 
\bP_{0,n}(1)^{[d]}\bP_{m,k}(a)=\sum_{j\le i}(-1)^{i-j}\frac{i!}{j!}{m\choose i}A_j(i,n)
\bP_{m-i,k+nj-i}(a{\bf a}_0^{i-j})\bP_{0,n}(1)^{[d-j]}).$$
\end{defi}

\begin{thm}\label{divided-powers-thm} The action of $\DD(\CH^*(\CC),K)$ on $\CH^*(\CC^{[\bullet]})$
extends to the action of $\wt{U}_1(\CH^*(\CC),K)$ (resp., $\wt{U}_2(\CH^*(\CC),K)$), such that
the action of $\bP_{n,0}(\CC)^{[d]}$ (resp., $\bP_{0,n}(\CC)^{[d]}$)
is given by \eqref{div-eq1} (resp., \eqref{div-eq2}). Furthermore, these relations hold on the level of correspondences. Also, similar statements hold for cohomology.
\end{thm}

First, we need to calculate some intersection-products.

\begin{lem}\label{inter-mult-lem2}
Recall that for each $m\ge 0$ we denote by $Z_{m,N}\sub \CC\times_S \CC^{[N]}$ the image of
$t_{m,N}:\CC\times_S\CC^{[N-m]}\to\CC\times_S\CC^{[N]}$.
For $m\le n$ we can consider the fine intersection-product
$[Z_{m,N}]\cdot [Z_{n,N}]\in\CH^m(Z_{n,N})$. We have the following formula:
$$[Z_{m,N}]\cdot [Z_{n,N}]=\sum_{i\ge 0}(-1)^i i!{m\choose i}{n\choose i} p_1^*K^i\cdot [Z_{m+n-i,N}],$$
where $p_1:Z_{n,N}\to \CC$ is the natural projection. 
\end{lem}

\Pf . We have an isomorphism $t_{n,N}:\CC\times \CC^{[N-n]}\wt{\to} Z_{n,N}$. 
Under this isomorphism the intersection-product in question becomes $t_{n,N}^*[Z_{m,N}]$, and the
required formula is equivalent to
$$t_{n,N}^*[Z_{m,N}]=\sum_{i\ge 0}(-1)^i{m\choose i}{n\choose i} p_1^*K^i\cdot [Z_{m-i,N-n}].$$
For $m=1$ this boils down to the identity
$$t_{n,N}^*[\DD_N]=\DD_{N-n}-n p_1^*K$$
that follows easily from \eqref{D-res-eq}. To deduce the case of general $m$ we use Lemma \ref{Z-lem}.
\ed

\begin{lem}\label{Pi-d-lem} 
Consider the cartesian diagram 
\begin{diagram}
\Pi_d(i_1,\ldots,i_n) &\rTo{} & \CC^{[d]}\\
\dTo{} & &\dTo{\De_n}\\
\CC\times_S \CC^{[d-i_1]}\times_S\ldots\times_S \CC^{[d-i_n]} &\rTo{s_{i_1,\ldots,i_n;d}}& (\CC^{[d]})^n,
\end{diagram}
where $\De_n$ is the diagonal embedding, while $s_{i_1,\ldots,i_n;d}$
is given by $(p,D_1,\ldots,D_n)\mapsto (i_1p+D_1,\ldots,i_np+D_n)$.
For each $j\ge\max(i_1,\ldots,i_n)$ we have a natural map
$$q^j:\CC\times_S \CC^{[d-j]}\to\Pi_d(i_1,\ldots,i_n):(p,D)\mapsto (p,(j-i_1)p+D,\ldots,(j-i_n)p+D,jp+D),$$
where we view $\Pi_d(i_1,\ldots,i_n)$ as a subvariety of
$\CC\times_S \CC^{[d-i_1]}\times_S\ldots\times_S \CC^{[d-i_n]}\times_S \CC^{[d]}$.
Set $i=\sum_{s=1}^n i_s$. Then we have
$$s_{i_1,\ldots,i_n;d}^![\CC^{[d]}]=
\sum_{j\ge 0} (-1)^j a(i_1,\ldots,i_n;j)\cdot p_1^*K^j\cdot q^{i-j}_*[C\times C^{[d-i+j]}]$$
where the coefficients $a(i_1,\ldots,i_n;j)$ are defined recursively by
\begin{equation}\label{a-rec-eq}
\begin{array}{l}
a(i_1,\ldots,i_n;j)=\sum_{k=0}^j k!{i_1\choose k}{i_2+\ldots+i_n-j+k\choose k}a(i_2,\ldots,i_n,j-k),\\
a(i_1;j)=\delta_{j,0}
\end{array}
\end{equation}
(note that $a(i_1,\ldots,i_n;j)=0$ unless $j\le i-\max(i_1,\ldots,i_n)$).
\end{lem}

\Pf . Note that for $n=1$ we have $\Pi_d(i_1)=\CC\times_S \CC^{[d-i_1]}$, and
the formula holds trivially. For $n>1$ we have the following commutative diagram with
cartesian squares: 
\begin{diagram}
\Pi_d(i_1,\ldots,i_n) &\rTo{} & \Pi_d(i_2,\ldots,i_n) &\rTo{} & \CC^{[d]}\\
\dTo{} & & \dTo{} & &\dTo{\De_n}\\
\CC\times_S \CC^{[d-i_1]}\times_S\ldots\times_S \CC^{[d-i_n]} 
&\rTo{t_{i_1,d}\times\id\times\ldots\times\id} & 
\CC\times_S \CC^{[d]}\times_S \CC^{[d-i_2]}\times_S\ldots\times_S \CC^{[d-i_n]} &\rTo{}& 
(\CC^{[d]})^n\\
\dTo{p_{12}} & &\dTo{p_{12}}\\
\CC\times_S \CC^{[d-i_1]} &\rTo{t_{i_1,d}} & \CC\times_S \CC^{[d]}
\end{diagram}
where the second arrow in the second row is induced by $s_{i_2,\ldots,i_n;d}$.
It follows that 
$$s_{i_1,\ldots,i_n;d}^![C^{[d]}]=t_{i_1,d}^!s_{i_2,\ldots,i_n;d}^![C^{[d]}].$$
Now the required equality and the recursive formula for the coefficients $a(i_1,\ldots,i_n;j)$ follow easily by induction using Lemma \ref{inter-mult-lem2}.
\ed

\begin{lem}\label{inter-mult-lem3}
Consider the cartesian square
\begin{diagram}
Z^{d,[n]}_{k,N} &\rTo{} & \CC\times_S \CC^{[N-k]}\\
\dTo{} & & \dTo{s_{k,N}}\\
\CC^{[d]}\times_S \CC^{[N-nd]} &\rTo{s^{[n]}_{d,N}} & \CC^{[N]}
\end{diagram}
where $s^{[n]}_{d,N}(D_1,D_2)=nD_1+D_2$.
For every $i,j$ such that $i\le k$ and $i\le nj$ we have a closed embedding
$$q^{i,j}:\CC\times_S \CC^{[d-j]}\times_S \CC^{[N-nd-k+i]}\to Z^{d,[n]}_{k,N}: (p,D_1,D_2)
\mapsto (jp+D_1,(k-i)p+D_2,p, (nj-i)p+nD_1+D_2),$$
where we view $Z^{d,[n]}_{k,N}$ as a subset of 
$\CC^{[d]}\times_S \CC^{[N-nd]}\times_S \CC\times_S \CC^{[N-k]}$.
Then we have the following formula for the intersection-product in the above diagram:
$$[\CC^{[d]}\times_S \CC^{[N-nd]}]\cdot [\CC\times_S \CC^{[N-k]}]=
\sum_{0\le j\le i\le k, i\le nj}(-1)^{i-j}\frac{i!}{j!}{k\choose i} A_j(i,n) 
q^{i,j}_*(p_1^*K^{i-j}\cdot [\CC\times_S \CC^{[d-j]}\times_S \CC^{[N-nd-k+i]}]),$$
where we use the numbers $(A_j(i,n))$ introduced in Lemma \ref{div-sum-lem}(ii).
\end{lem}

\Pf . 
Consider the following commutative diagram with
cartesian squares: 
\begin{diagram}
Z^{d,[n]}_{k,N} &\rTo{} & \Pi_{k,(d)^n,N} &\rTo{} & \CC\times_S \CC^{[N-k]}\\
\dTo{} & & \dTo{} & &\dTo{}\\
\CC^{[d]}\times_S \CC^{[N-nd]} &\rTo{\De_n\times\id} & 
(\CC^{[d]})^n\times_S \CC^{[N-nd]} &\rTo{\a}& \CC^{[N]}\\
\dTo{} & &\dTo{}\\
\CC^{[d]} &\rTo{\De_n}&(\CC^{[d]})^n
\end{diagram}
where $\De_n$ is the diagonal embedding, the map $\a$
is given by the addition of divisors. We have to calculate
$\De_n^!\a^![\CC\times_S \CC^{[N-k]}]$.
Iterating Lemma \ref{mult-lem2} we obtain 
$$\a^![\CC\times_S \CC^{[N-k]}]=\sum_{i_1+\ldots+i_n=i\le k}\frac{k!}{i_1!\ldots i_n!(k-i)!}
[\CC\times_S \CC^{[d-i_1]}\times_S\ldots\times_S \CC^{[d-i_n]}\times_S \CC^{[N-nd-k+i]}].$$
Next, by Lemma \ref{Pi-d-lem}
\begin{align*}
&\De_n^![\CC\times_S\CC^{[d-i_1]}\times_S\ldots\times_S \CC^{[d-i_n]}]=
s_{i_1,\ldots,i_n;d}^![\CC^{[d]}]=\\
&\sum_{l\ge 0}(-1)^l a(i_1,\ldots,i_n;l) p_1^*K^l\cdot [\CC\times_S\CC^{[d-i+l]}],
\end{align*}
where $i=i_1+\ldots+i_n$.
It follows that
$$\De_n^!\a^![\CC\times_S \CC^{[N-k]}]=\sum_{0\le l\le i\le k}(-1)^l{k\choose i}
b(i,l;n) p_1^*K^l\cdot [\CC\times_S\CC^{[d-i+l]}\times_S\CC^{[N-nd-k+i]}],$$
where 
$$b(i,l;n)=\sum_{i_1+\ldots+i_n=i}\frac{i!}{i_1!\ldots i_n!}a(i_1,\ldots,i_n;l).$$
It remains to show that $b(i,l;n)=\frac{i!}{(i-l)!}A_{i-l}(i,n)$. To this end we use the recursive formulas
$$b(i,l;n)=\sum_{0\le i_1\le i,0\le k\le l}\frac{i!(i-i_1-l+k)!}{(i-i_1)!(i_1-k)!(i-i_1-l)!k!}b(i-i_1,l-k;n-1),
\ \ \ b(i,l;1)=\de_{l,0}$$
that follow immediately from \eqref{a-rec-eq}. Note that from the formula 
${n\choose i}={n-1\choose i}+{n-1\choose i-1}$ one can derive a similar recursive formula for 
$A_d(i,n)$:
$$A_d(i,n)=\sum_{r+s+t=d}\frac{d!}{r!s!t!}A_{d-r}(i-r-s,n-1).$$
Now the required equality follows easily by induction in $n$.
\ed

\noindent
{\it Proof of Theorem \ref{divided-powers-thm}.}
Relations of type (ii) are easy to check, so we will concentrate on relations of type (iii).
In the notation of Lemma \ref{inter-mult-lem3}
the composition $P_{m,k}(a)\circ P_{n,0}(\CC)^{[d]}$ acting on $\CC^{[N-nd]}$ is given by
$x\mapsto q_{2*}(w\cdot p_{\CC}^*a\cdot q_1^*x)$, where 
$w\in\CH^*(Z^{d,[n]}_{k,N})$ is the intersection-product computed in this lemma,
$q_1$ is the composition 
$$Z^{d,[n]}_{k,N}\to \CC^{[d]}\times_{\SS} \CC^{[N-nd]}\stackrel{p_2}{\to}\CC^{[N-nd]},$$
$q_2$ is the composition
$$
\begin{diagram}
Z^{d,[n]}_{k,N} &\rTo &\CC\times_S\CC^{[N-k]} &\rTo{s_{m,N-k+m}} &\CC^{[N-k+m]},
\end{diagram}
$$
and $p_{\CC}:Z^{d,[n]}_{k,N}\to\CC$ is the natural projection.
The formula for $w$ leads to the expression for the $P_{m,k}(a)\circ P_{n,0}(\CC)^{[d]}$ as the
linear combination of the operators defined by the cycles $p_{\CC}^*(a\cdot K^{i-j})$ over the correspondences
$$
\begin{diagram}
& & \CC\times_S\CC^{[d-j]}\times_S\CC^{[N-nd-k+i]} & \\
&\ldTo{s_{k-i,N-nd}p_{13}} & &\rdTo{q} &\\
\CC^{[N-nd]} & & & &\CC^{[N-k+m]}
\end{diagram}
$$
where $q(p,D_1,D_2)=(m+nj-i)p+nD_1+D_2$. It is easy to see that the same correspondence
arises when computing $P_{n,0}(\CC)^{[d-j]}\circ P_{m+nj-i,k-i}(a\cdot K^{i-j})$ and that the coefficients
match. The second relation of type (iii) is checked similarly: the operators involved in it are defined
by the transposes of the above correspondences.
\ed

\begin{ex} In the case of a trivial family $\CC=C\times S$ the relations established in 
Theorem \ref{divided-powers-thm} take form
$$P_{m,k}(a)P_{n,0}(C\times S)^{[d]}=
\sum_{i=0}^d{k\choose i} n^i P_{n,0}(C\times S)^{[d-i]}P_{m+i(n-1),k-i}(a),$$
$$P_{0,n}(C\times S)^{[d]}P_{m,k}(a)=\sum_{i=0}^d{m\choose i}n^i P_{m-i,k+i(n-1)}(a)
P_{0,n}(C\times S)^{[d-i]},$$
where $a\in\CH^*(C\times S)$.
\end{ex}

\begin{rems} 1. One should be able to establish also some commutation relations between
$P_{n_1,0}(\CC)^{[d_1]}$ and $P_{0,n_2}(\CC)^{[d_2]}$. The simplest example of such relations
is given in Proposition \ref{div-com-prop} below.

\noindent 2. For other operators $P_{k,m}(\CC)$ 
one should be able construct some modified divided powers.
For example, we have $P_{1,1}(\CC)(x)=Nx$
for $x\in\CH^*(\CC^{[N]})$, so one cannot construct $P_{1,1}(\CC)^2/2$, however, one can construct
$(P_{1,1}(\CC)^2-P_{1,1}(\CC))/2$.
\end{rems}

\begin{prop}\label{div-com-prop}
All the operators in the family $\{P_{1,0}(\CC)^{[d]}\ |\ d\ge 1\}\cup
\{P_{0,1}(\CC)^{[d]}\ |\ d\ge 1\}$ commute with each other (on the level of correspondences).
\end{prop}

\Pf . The fact that the operators within each set commute with each other follows from the commutativity
of the Pontryagin product. Now let us check that $P_{1,0}(\CC)^{[d_1]}$ commutes with
$P_{0,1}(\CC)^{[d_2]}$.
The composition $P_{0,1}(\CC)^{[d_2]}\circ P_{1,0}(\CC)^{[d_1]}$ acting on $\CH^*(\CC^{[N]})$
is given by the following correspondence $\Pi$ from $\CC^{[N]}$ to $\CC^{[N+d_1-d_2]}$ equipped
with a class of dimension $N+d_1$:
$$\Pi=\{(D_1,E_1,D_2,E_2)\in\CC^{[d_1]}\times_S\CC^{[N]}\times_S\CC^{[d_2]}\times_S
\CC^{[N+d_1-d_2]}\ |\ D_1+E_1=D_2+E_2\},$$
where the map $q:\Pi\to\CC^{[N]}\times_S\CC^{[N+d_1-d_2]}$ sends $(D_1,E_1,D_2,E_2)$ to
$(E_1,E_2)$. This correspondence is equipped with the natural intersection-product class of dimension
$N+d_1$
\begin{equation}\label{div-com-inter-class}
[\CC^{[d_1]}\times_S\CC^{[N]}]\cdot [\CC^{[d_2]}\times_S\CC^{[N+d_1-d_2]}].
\end{equation}
On the other hand, the composition $P_{1,0}(\CC)^{[d_1]}\circ P_{0,1}(\CC)^{[d_2]}$ is given by the correspondence
$$
\begin{diagram}
& &\CC^{[d_2]}\times_S\CC^{[N-d_2]}\times_S\CC^{[d_1]}& &\\
&\ldTo{\a_{d_2,N-d_2}p_{12}}& &\rdTo{\a_{N-d_2,d_1}p_{23}}&\\
\CC^{[N]}& & & &\CC^{[N+d_1-d_2]}
\end{diagram}
$$
The natural map 
$$\CC^{[d_2]}\times_S\CC^{[N-d_2]}\times_S\CC^{[d_1]}\to\Pi:(D_2,E',D_1)\mapsto
(D_1,D_2+E',D_2,D_1+E')$$
is an isomorphism onto an irreducible component $\Pi_0\sub\Pi$.
Other irreducible components $\Pi_i\sub\Pi$ are numbered by $i$, such that $1\le i\le\min(d_1,d_2)$.
Namely, $\Pi_i$ is the image of the map
$$\CC^{[i]}\times_S\CC^{[d_2-i]}\times_S\CC^{[N-d_2+i]}\times_S\CC^{[d_1-i]}\to\Pi:
(D_0,D'_2,E',D'_1)\mapsto (D_0+D'_1,D'_2+E', D_0+D'_2,D'_1+E').$$
Since the composition of this map with $q:\Pi\to\CC^{[N]}\times_S\CC^{[N+d_1-d_2]}$
factors through $\CC^{[d_2-i]}\times_S\CC^{[N-d_2+i]}\times_S\CC^{[d_1-i]}$, we derive
that $q(\Pi_i)$ has dimension $\le N+d_1-i$. Therefore, the components $\Pi_i$ with $i\ge 1$
give zero contribution to the composition $P_{0,1}(\CC)^{[d_2]}\circ P_{1,0}(\CC)^{[d_1]}$.
It remains to prove that $[\Pi_0]$ appears in the intersection-product \eqref{div-com-inter-class}
with multiplicity $1$. To this end we can replace $\CC^{[d_2]}$ by the cartesian product $\CC^{d_2}$
and then use iteratively Lemma \ref{mult-lem2} (in the case $m=1$).
\ed

Let us denote by $\DD_{t,u,\Q}=\Q[t,u,\pa_t,\pa_u]$ the algebra of differential operators in two variables.
Let us also denote by 
$$\Z[t,u^{[\bullet]}, \pa_t^{[\bullet]}, \pa_u]\sub\DD_{t,u,\Q}$$ 
the subalgebra over $\Z$ generated by $t$, $\pa_u$ and by the divided powers of $u$ and $\pa_t$. 
We will also consider the $\Z$-subalgebras $\Z[t,\pa_t^{[\bullet]}]$ and $\Z[u^{[\bullet]},\pa_u]$ in
this algebra.

\begin{cor}\label{Heis-cor} 
Assume that we are given a point $p_0\in\CC(S)$. Then
there is an action of $\Z[t,u^{[\bullet]}, \pa_t^{[\bullet]}, \pa_u]$
on $\CH^*(\CC^{[\bullet]})$ such that
\begin{align*}
&t\mapsto P_{1,0}([p_0(S)])+\psi\cdot P_{1,0}(\CC) \ \ 
(\text{Pontryagin product with }[p_0(S)]+\psi\cdot[\CC])
\\
&u^{[d]}\mapsto P_{1,0}(\CC)^{[d]} \ \ (\text{Pontryagin product with } [\CC^{[d]}]),\\ 
&\pa_t^{[d]}\mapsto P_{0,1}(\CC)^{[d]},\\
&\pa_u\mapsto P_{0,1}([p_0(S)]), 
\end{align*}
where $\psi=p_0^*K\in\CH^1(S)$.
\end{cor}

Here are some simple observations on actions of 
$\Z[t,u^{[\bullet]}, \pa_t^{[\bullet]}, \pa_u]$ (probably well known).

\begin{lem}\label{div-mod-lem} 
Let $M$ be a $\Z[u^{[\bullet]},\pa_u]$-module such that for every $x\in M$ one has
$\pa_u^nx=0$ for all $n\gg 0$ (resp., $\Z[t,\pa_t^{[\bullet]}]$-module such that for every $x\in M$ one has
$\pa_t^{[n]}x=0$ for all $n\gg 0$).
Set $M_0=\{x\in M\ |\ \pa_ux=0\}$ (resp., $M_0=\{x\in M\ |\ \pa_t^{[n]}x=0 \text{ for all }n>0\}$).
Then the submodule $M_0[u^{[\bullet]}]\sub M$ (resp., $M_0[t]\sub M$) consisting of elements
of the form $\sum a_i u^{[i]}$ (resp., $\sum a_i t^i$) with $a_i\in M_0$, coincides with the entire $M$.
\end{lem}

\Pf . Assume first that $M$ is a module over $\Z[u^{[\bullet]},\pa_u]$.
For $x\in M$ let $n$ be the minimal number such that $\pa_u^n(x)\in M_0[u^{[\bullet]}]$.
Assume that $n>0$. Let
$$\pa_u^n(x)=a_0+a_1u+\ldots+a_ku^{[k]}.$$
where $a_i\in M_0$. Then 
$$\pa_u(\pa_u^{n-1}(x)-a_0u-a_1u^{[2]}-\ldots-a_ku^{[k+1]})=0,$$
i.e., $\pa_u^{n-1}(x)-a_0u-a_1u^{[2]}-\ldots-a_ku^{[k+1]}\in M_0$. It follows that
$\pa_u^{n-1}(x)\in M_0[u^{[\bullet]}]$ contradicting the choice of $n$.

Now, let $M$ be a module over $\Z[t,\pa_t^{[\bullet]}]$. For $x\in M$ let $n$ be the minimal number
such that $\pa_t^{[k]}(x)\in M_0[t]$ for all $k\ge n$. 
Assume that $n>0$ and let us lead this to contradiction.
Replacing $x$ by $\pa_t^{[n-1]}(x)$ we can reduce the proof to the case $n=1$.
In this case set for $k\ge 1$
$$f_k=\pa_t^{[k]}(x)=\sum_{i\ge 0}a_{k,i}t^i\in M_0[t].$$
Note that $f_k=0$ for $k\gg 0$, so we can form the finite sum
$$f=\sum_{k\ge 1}a_{k,0}t^k\in M_0[t].$$
Using the identities $\pa_t^{[m-k]}f_k={m\choose k}f_m$ one can easily check that
$\pa_t^{[k]}f=f_k$ for every $k\ge 1$. Hence, $x-f\in M_0$, i.e., $x\in M_0[t]$.
\ed

For every abelian group $\Ga$ there is a natural structure of 
a $\Z[t,u^{[\bullet]}, \pa_t^{[\bullet]}, \pa_u]$-module on
$\Ga[t,u^{[\bullet]}]$ such that $\pa_t^{[n]}(\ga t^iu^{[j]})=\ga{i\choose n}t^{i-n}u^{[j]}$. 
Note that the operator $\pa_u$ on this module is surjective.

\begin{prop}\label{Heis-mod-prop} 
Let $M$ be a $\Z[t,u^{[\bullet]}, \pa_t^{[\bullet]}, \pa_u]$-module such that for every $x\in M$ one has 
$\pa_t^{[n]}x=\pa_u^nx=0$ for all $n\gg 0$. Then 

\noindent (i) $M\simeq M_0[t,u^{[\bullet]}]$ as $\Z[t,u^{[\bullet]}, \pa_t^{[\bullet]}, \pa_u]$-module, where 
$$M_0=\{x\in M\ |\ \pa_ux=0, \pa_t^{[n]}x=0 \text{ for all }n>0\};$$

\noindent (ii) the operator $t:M\to M$ is injective, and the operator $\pa_u:M\to M$ is surjective.
\end{prop}

\Pf . (i) First, it is easy to check that any submodule of an $\Z[t,u^{[\bullet]}, \pa_t^{[\bullet]}, \pa_u]$-module of the form $\Ga[t,u^{[\bullet]}]$ itself
has the form $\Ga'[t,u^{[\bullet]}]$ for a subgroup $\Ga'\sub\Ga$. Indeed, this follows easily
from the fact that $\pa_t^{[m]}\pa_u^n(\ga t^m u^{[n]})=\ga$ for $\ga\in\Ga$.
Therefore, the natural morphism of $\Z[t,u^{[\bullet]}, \pa_t^{[\bullet]}, \pa_u]$-modules 
$M_0[t,u^{[\bullet]}]\to M$ is injective (since $M_0$ embeds into $M$).
It remains to prove that this morphism is also surjective. Viewing $M$ as a
$\Z[u^{[\bullet]},\pa_u]$-module and applying Lemma \ref{div-mod-lem}, we derive that
$M=(\ker \pa_u)[u^{[\bullet]}]$. Next, applying the same lemma to the $\Z[t,\pa_t^{[\bullet]}]$-module
$\ker\pa_u$, we obtain $\ker\pa_u=M_0[t]$. Hence, $M=M_0[t,u^{[\bullet]}]$.

\noindent
(ii) This follows from (i) since for a module of the form $\Ga[t,u^{[\bullet]}]$ 
injectivity of $t$ (resp., surjectivity of $\pa_u$) is clear.
\ed

In the case when $S$ is a point the following result is due to Collino~\cite{Col2}.

\begin{prop}\label{Col-prop} Let $i_{N}:\CC^{[N-1]}\to\CC^{[N]}$ be the closed embedding associated
with a point $p_0\in\CC(S)$. Then the homomorphism $i_{N*}:\CH^*(\CC^{[N-1]})\to\CH^*(\CC^{[N]})$
(resp., $i_N^*:\CH^*(\CC^{[N]})\to\CH^*(\CC^{[N-1]})$) is injective (resp., surjective).
\end{prop}

\Pf . Note that $i_{N*}=P_{1,0}([p_0(S)])$ and $i_N^*=P_{0,1}([p_0(S)])$. Now the
surjectivity of $P_{0,1}([p_0(S)])$ follows immediately from Corollary \ref{Heis-cor} and Proposition 
\ref{Heis-mod-prop}(ii). To deal with injectivity of $P_{1,0}([p_0(S)])$ we 
can modify the action of Corollary \ref{Heis-cor}(i) as follows:
\begin{align*}
&t\mapsto P_{1,0}(p_0(S)),\\
&u^{[d]}\mapsto P_{1,0}(\CC)^{[d]},\\ 
&\pa_t^{[d]}\mapsto P_{0,1}(\CC)^{[d]},\\
&\pa_u\mapsto P_{0,1}([p_0(S)])+\psi\cdot P_{0,1}(\CC).
\end{align*}
It remains to apply Proposition \ref{Heis-mod-prop}(ii) to this action.
\ed

We also get the following corollary from Proposition \ref{Heis-mod-prop}.

\begin{cor}\label{module-cor}
For a point $p_0\in \CC(S)$ consider the action of $\Z[t,u^{[\bullet]}, \pa_t^{[\bullet]}, \pa_u]$
on $\CH^*(\CC^{[\bullet]})$ given by Corollary \ref{Heis-cor}. Then
there is an isomorphism of $\Z[t,u^{[\bullet]}, \pa_t^{[\bullet]}, \pa_u]$-module
$$\CH^*(\CC^{[\bullet]})\simeq K_{p_0}[t,u^{[\bullet]}],$$
where 
$$K_{p_0}=\{x\in\CH^*(\CC^{[\bullet]})\ |\ P_{0,1}(p_0(S))x=0, 
P_{0,1}(\CC))^{[d]}x=0 \text{ for all }d\ge 1\}.$$
\end{cor}

\begin{rems}
1. The isomorphism of the above corollary is compatible with the bigrading of $\CH^*(C^{[\bullet]})$:
$t$ (resp., $u$) sends $\CH^i(C^{[N]})$ to $\CH^{i+1}(C^{[N+1]})$ (resp., $\CH^i(C^{[N+1]})$).
In the case $S=\Spec(k)$ the strong stability conjecture (see \cite{KV}, 2.13) is equivalent to the condition that $K_{p_0}\cap\CH^p(C^{[N]})$ is a torsion group for $N>2p$.

\noindent 2. In the case $S=\Spec(k)$ the decomposition of $\CH^*(C^{[\bullet]})$ into the
direct summands of the form $K_{p_0}t^mu^{[n]}$ is the $\Z$-version of the well known motivic decomposition over $\Q$ obtained by using $\lambda$-operations (see \cite{dB1}).

\noindent 3. Using the divided powers of $P_{0,1}(C)$ we can avoid tensoring with $\Q$ in the
proof of Theorem \ref{curve-thm}. Instead one has to use the fact that $P_{0,1}(C)^{[d]}$ 
commutes with $P_{1,0}(a)$ for $a\in A_0(C)$, and hence the subgroup
$A_0(C)^{*n}\sub\CH_0(C^{[n]})$ is killed by all the operators $P_{0,1}(C)^{[d]}$.
\end{rems}

Using Corollary \ref{Heis-cor} we get an interesting $\splin_2$-action on the motive of $\CC^{[N]}$.
In the case when $S$ is a point we obtain in this way a Lefschetz $\splin_2$-action on 
$\CH^*(C^{[N]})$.

\begin{thm}\label{Lefschetz-thm} 
(i) Fix a point $p_0\in \CC(S)$. Then for every $N\ge 0$ the operators
\begin{align*}
&e(x)=[\RR]\cdot x+\psi\cdot P_{1,0}(\CC)P_{0,1}([p_0(S)])(x),\\
&f=P_{1,0}(\CC)P_{0,1}(\CC),\\
&h=P_{1,0}([p_0(S)])P_{0,1}(\CC)-P_{1,0}(\CC)P_{0,1}([p_0(S)])+\psi\cdot P_{1,0}(\CC)P_{0,1}(\CC)
\end{align*}
(given by algebraic correspondences)
define compatible actions of the Lie algebra $\splin_2$ on $\CH^*(\CC^{[N]})$ and on 
$H^*(\CC^{[N]},\Q)$,
where $\RR=\RR_N\sub\CC^{[N]}$ is the divisor associated with $p_0$ (see section \ref{Jac-sec}). 

\noindent
(ii) In the case when $S$ is a point (so we write $\CC=C$)
the operator $h$ acts as $(i-N)\id$ on $H^i(C^{[N]},\Q)$. In other words, the
action of $(e,f,h)$ on $H^*(C^{[N]},\Q)$ is the Lefschetz action corresponding to the ample divisor
$R\sub C^{[N]}$ (associated with $p_0$).
\end{thm}

\Pf . (i) Consider the action of $\Z[t,u^{[\bullet]}, \pa_t^{[\bullet]}, \pa_u]$ on $\CH^*(\CC^{[\bullet]})$
given by Corollary \ref{Heis-cor}. Then the operators $e=t\pa_u$, $f=u\pa_t$ and $h=t\pa_t-u\pa_u$
satisfy the relations of $\splin_2$.
By definition, we have 
$$P_{1,0}([p_0])|_{\CC^{[N-1]}}=i_{N*}, \ \ P_{0,1}([p_0])|_{\CC^{[N]}}=i_{N}^*,$$ 
where $i_N:\CC^{[N-1]}\to\CC^{[N]}$ is
the embedding associated with $p_0$.
Hence, for $x\in\CH^*(\CC^{[N]})$ one has
$$P_{1,0}([p_0])P_{0,1}([p_0])(x)=i_{N*}i_N^*x=[\RR]\cdot x.$$
This implies our formula for $e$.

\noindent
(ii) It is easy to see that two $\splin_2$-triples $(e,f,h)$ and $(e',f',h')$ acting
on the same finite-dimensional space $V$ such that $e=e'$ and $[h,h']=0$ necessarily coincide
(i.e., $h=h'$ and $f=f'$). Indeed, consider the decomposition $V=\bigoplus_{m,n} V_{m,n}$,
where $h$ (resp., $h'$) acts by $m\cdot\id$ (resp., $n\cdot\id$) on $V_{m,n}$.
Consider also the subspace $W=\ker(e)\sub V$ and the induced decomposition 
$W=\bigoplus_{m\ge 0,n\ge 0} W_{m,n}$, where $W_{m,n}=W\cap V_{m,n}$.
Viewing $V$ as a representation of $(e,f,h)$ we deduce that 
$$e^iV\cap W=\bigoplus_{m\ge i,n}W_{m,n}\text{ for all }i\ge 0.$$
On the other hand, using the triple $(e=e',f',h')$ we get 
$$e^iV\cap W=\bigoplus_{m,n\ge i}W_{m,n}\text{ for all }i\ge 0.$$
This immediately implies that $W_{m,n}=0$ for $m\neq n$, i.e., $h=h'$. It is well known that
this implies $f=f'$.

Now comparing the action of our operators $(e,f,h)$ on $H^*(C^{[N]},\Q)$ with the Lefschetz action corresponding to $[R]$ and using the above observation we conclude that these two actions coincide.
\ed

\begin{rem} It is well known that the standard conjecture B of Lefschetz type for $X$ and 
for a projective bundle over $X$ are equivalent (see \cite{Li}).
Also, if it is true for some variety $X$ then
it is true for an ample divisor in $X$. Thus, standard conjecture B for all $C^{[N]}$ and for the Jacobian
$J$ are equivalent.
In particular, from the above theorem we get a new proof of the standard conjecture B for $J$. 
\end{rem}

Using our methods we easily recover the following result proved by S.~del Ba\~{n}o in
\cite{dB2}.

\begin{cor}\label{Lefschetz-cor} 
The hard Lefschetz theorem for the operator $L$ of multiplication by the class
$[R]\in H^2(C^{[N]},\Z)$ (associated with a point $p_0\in C$) holds over $\Z$, i.e., the map
$$L^i:H^{N-i}(C^{[N]},\Z)\to H^{N+i}(C^{[N]},\Z)$$
is an isomorphism for all $i\ge 0$.
\end{cor}

\Pf . Since the action of $\Z[t,u^{[\bullet]}, \pa_t^{[\bullet]}, \pa_u]$ defined in Corollary \ref{Heis-cor} is given by algebraic correspondences, it can also be defined
for $H^*(C^{[N]},\Z)$. The operator $L$ corresponds to the action of $e=t\pa_u$. As we have seen
in Theorem \ref{Lefschetz-thm}, the operator $h=t\pa_t-u\pa_u$ acts as $(i-N)\id$ on $H^i(C^{[N]},\Z)$.
On the other hand, by Proposition \ref{Heis-mod-prop}(i) we have an isomorphism of $H^*(C^{[N]},\Z)$
with a $\Z[t,u^{[\bullet]}, \pa_t^{[\bullet]}, \pa_u]$-module of the form 
$\Ga[t,u^{[\bullet]}]$. For such a module we have
$h(\ga t^mu^{[n]})=(m-n)\ga t^mu^{[n]}$, and our assertion follows from the formula
$$e^{n-m}(\ga t^mu^{[n]})=\ga t^nu^{[m]} \text{ for }n\ge m.$$
\ed

 \section{Tautological cycles}
\label{taut-sec}

Similarly to the case of cycles on the Jacobian considered in \cite{Bmain} and \cite{P-lie}, 
we are going to 
define the subalgebra of tautological classes in $\CH^*(\CC^{[\bullet]})$ (resp., $\CH^*(\JJ)$).
We will show that the operators $(P_{m,k}(a))$ (resp., $T_k(m,a)$) act on this subalgebra by some differential operators. 

We start with an abstract setup.
Let $R$ be a supercommutative ring,
$A$ a supercommutative $R$-algebra with a fixed even element ${\bf a}_0\in A$,
$\DD(A,{\bf a}_0)$ the corresponding Lie superalgebra (see the Introduction). 
Let us also denote by $\DD_+(A,{\bf a}_0)\sub\DD(A,{\bf a}_0)$ the subalgebra generated by the operators
$\bP_{m,k}(a)$ with $m\ge k$.

\begin{prop}\label{diff-op-prop} Let $B$ be  a supercommutative $R$-algebra,
and let $x_m:A\to B$, $m\ge 0$, be a family of even $R$-linear maps. 

\noindent
(i) Assume that $\DD(A,{\bf a}_0)$ (resp., $\DD_+(A,{\bf a}_0)$)
acts $R$-linearly on $B$ in such a way that 
$$\bP_{m,0}(a)(b)=x_m(a)\cdot b$$
for all $m\ge 0$ and $a\in A$.
Assume also that 
$$\bP_{m,k}(a)(1)=0 \text{  for }k>0$$
(resp., for $m\ge k>0$).
Let us denote by $\TT B\sub B$ 
the $R$-subalgebra generated by $(x_m(a))$. 
Then $\TT B$ is stable under the action of $\DD(A,{\bf a}_0)$ (resp., $\DD_+(A,{\bf a}_0)$). 
Furthermore, if we view $\TT B$ as the quotient of the superalgebra of polynomials 
$R[x_m(a)\ |\ m\ge 0, a\in G]$,
where $G$ is some set of homogeneous generators of $A$ as an $R$-module,
then the action of this Lie algebra on $\TT B$ is given by the formulas
\begin{align*}
&\bP_{m,k}(a)=\sum_{s\ge 0;k_1+\ldots+k_s=k, k_i\ge 1;
n_1,\ldots,n_s\ge 0; a_1,\ldots,a_s\in G}
(-1)^{k-s}\frac{k!}{s!}{n_1\choose k_1}\cdot\ldots\cdot{n_s\choose k_s}\times\\
&x_{m-k+n_1+\ldots+n_s}(a a_s\ldots a_1\cdot {\bf a}_0^{k-s})\pa_{x_{n_1}(a_1)}\ldots\pa_{x_{n_s}(a_s)}, 
\end{align*}
where the case $s=0$ occurs only for $k=0$.

\noindent
(ii) Assume in addition that the above action of $\DD(A,{\bf a}_0)$ on $B$ extends to an action of
$\wt{U}_1(A,{\bf a}_0)$, so that 
$$\bP_{m,0}(1)^{[d]}(b)=x_m(1)^{[d]}\cdot b$$
for some elements $x_m(1)^{[d]}\in B$. Let $\wt{\TT} B$ be the subalgebra generated by $\TT B$ and
by all the elements $x_m(1)^{[d]}$. Then $\wt{\TT} B$ is stable under the action of $\wt{U}_1(A,{\bf a}_0)$,
and the action of $\bP_{m,k}(a)$ on it is given by the same formula as above (with $\pa_{x_m(1)}$ extended to the divided powers).
Similarly, if the action of $\DD(A,{\bf a}_0)$ on $B$ extends to an action of $\wt{U}_2(A,{\bf a}_0)$ 
such that $\bP_{0,n}(1)^{[d]}(1)=0$ for $n\ge 1, d\ge 1$, then
$\TT B$ is stable under the action of $\wt{U}_2(A,{\bf a}_0)$, and the action of
$\bP_{0,n}(1)^{[d]}$ on $\TT B$ is given by a differential operator of order $nd$.
\end{prop}

\Pf . (i) For $k=0$ we have $\bP_{m,0}(a)=x_m(a)$ by assumption.
The general case follows by induction in $k$ using the (super)commutator formula
$$[\bP_{m,k}(a),x_n(a')]=[\bP_{m,k}(a),\bP_{n,0}(a')]=
\sum_{i\ge 1}(-1)^{i-1}i!\cdot {k\choose i}{n\choose i}\bP_{m+n-i,k-i}(a\cdot a'\cdot {\bf a}_0^{i-1})$$
together with the assumption that $\bP_{m,k}(a)(1)=0$ for $k>0$.
Indeed, it is straightforward to check that the similar commutation relation 
holds for the differential operators in the right-hand side of the required formula.

\noindent (ii) Same proof as in (i) using the commutation relations in $\wt{U}_1(A,{\bf a}_0)$ (resp.,
$\wt{U}_2(A,{\bf a}_0)$).
\ed

Proposition \ref{diff-op-prop} can be applied to the algebra
$B=\CH_*(\CC^{[\bullet]})$ (equipped with Pontryagin product) and the operators 
$(P_{m,k}(a))$ acting on it. This leads to a definition of the subalgebra of tautological classes.

\begin{defi} The {\it (big) algebra of tautological classes on $\CC^{[\bullet]}$} 
$$\TCH^*(\CC^{[\bullet]})\sub\CH^*(\CC^{[\bullet]})$$
is the $\CH^*(S)$-subalgebra with respect to the Pontryagin product 
generated by all the classes of the form
$\De_{n*}(a)\in\CH^*(\CC^{[n]})$, where $\De_n:\CC\to\CC^{[n]}$ is the diagonal embedding,
$a\in\CH^*(\CC)$, $n\ge 1$. Let us also denote by 
$$\wt{\TCH}^*(\CC^{[\bullet]})\sub\CH^*(\CC^{[\bullet]})$$ 
the subalgebra generated by $\TCH^*(\CC^{[\bullet]})$ along with all the classes
$\de_n^{[d]}\in\CH^{(n-1)d}(\CC^{[nd]})$ (see section \ref{div-sec}).
Replacing Chow groups by cohomology we also define the
subalgebra of tautological classes $\TH^*(\CC^{[\bullet]})\sub H^*(\CC^{[\bullet]})$.
\end{defi}

We can also mimic the above definition of tautological classes in the case of the relative Jacobian
$\JJ$ for a family of curves $\CC/S$ equipped with a point $p_0\in\CC(S)$.

\begin{defi} The {\it (big) algebra of tautological classes on $\JJ$} 
$$\TCH^*(\JJ)\sub\CH^*(\JJ)$$
is the $\CH^*(S)$-subalgebra with respect to the Pontryagin product 
generated by all the classes of the form $[n]_*(\iota_*a)$, where 
$\iota=\si_1:\CC\to\JJ$ is the embedding
associated with $p_0$, $a\in\CH^*(\CC)$, $n\in\Z$.
Here we view $\CH^*(\JJ)$ as a $\CH^*(S)$-algebra via the homomorphism
$e_*:\CH^*(S)\to\CH^*(\JJ)$ associated with the neutral element $e\in\JJ(S)$.
Similarly, we define the subalgebra of tautological classes in cohomology.
\end{defi}

\begin{rem}
One can also consider smaller algebras of tautological classes by choosing a $\CH^*(S)$-subalgebra
$\AA\sub\CH^*(\CC)$ and considering only the classes $\De_{n*}(a)$ with $a\in\AA$.
For example, one can take as $\AA$ the subalgebra generated by $[p_0(S)]\in\CH^1(\CC)$,
or by the relative canonical class $K\in\CH^1(\CC)$, or by some other divisor class. In \cite{P-lie} we
worked with the subalgebra generated by $\chi+[p_0]$, where $2\chi=K$ (in the case $S=\Spec(k)$).
\end{rem}

\begin{thm}\label{taut-thm}
(i) The operators $(P_{m,k}(a))$ from the Introduction preserve the 
subalgebra of tautological classes $\TCH^*(\CC^{[\bullet]})\sub\CH^*(\CC^{[\bullet]})$
(resp., $\wt{\TCH}^*(\CC^{[\bullet]})$)
and act on it by the differential operators given in Proposition \ref{diff-op-prop} with
$x_n(a)=\De_{n*}(a)$. 

\noindent
(ii) The operators $T_k(m,a)$ from section \ref{Jac-sec}
preserve the subalgebra of tautological classes
$\TCH^*(\JJ)\sub\CH^*(\JJ)$ and act on it by differential operators (with respect to the Pontryagin product), so that $T_k(m,a)$ acts by a differential operator of order $k$. 

\noindent
(iii) The space of tautological classes with rational coefficients $\TCH^*(\JJ)_{\Q}\sub\CH^*(\JJ)_{\Q}$
is closed under the usual product and under the Fourier transform. It coincides with the 
$\CH^*(S)_{\Q}$-subalgebra
in $\CH^*(\JJ)_{\Q}$ with respect to the usual product, generated by the classes $\tau_k(a)$, $k\ge 0$,
$a\in\CH^*(\CC)$ (see \eqref{tau-eq}).

\noindent
(iv) For every $N$ the map $\si_{N*}:\CH^*(\CC^{[N]})\to\CH^*(\JJ)$
(resp., $\si_N^*:\CH^*(\JJ)_{\Q}\to\CH^*(\CC^{[N]})_{\Q}$) sends
tautological classes to tautological classes (resp., with rational coefficients).

\noindent
(v) Similar statements hold for the cohomology. 
\end{thm}

\Pf . (i) This follows immediately from Proposition \ref{diff-op-prop}.

\noindent (ii) The operator $T_0(m,a)$ is simply the Pontryagin product with 
$[m]_*\iota_*a\in\CH^*(\JJ)$.
The commutation relation of Theorem \ref{relations-thm} for $k'=0$ and $k\ge 1$ gives
\begin{align*}
&[T_k(m,a),T_0(m',a')]+\sum_{i\ge 1}\psi^i{k\choose i}m^{\prime i}T_{k-i}(m,a)T_0(m',a')=\\
&\sum_{i\ge 1}(-1)^{i-1}{k\choose i}m^{\prime i}T_{k-i}(m+m',aa'(K+2[p_0(S)])^{i-1})-
\psi^{k-1}m^{\prime k}p_0^*(a)T_0(m',a')\\
&-\sum_{i\ge 1}{k\choose i}m^{\prime i}\psi^{i-1}p_0^*(a')T_{k-i}(m,a).
\end{align*}
Since $T_k(m,a)(e_*x)=0$ this implies the assertion by induction in $k$.

\noindent
(iii) First of all, note that $g$-th Pontryagin power of $\si_{1*}(\CC)$ is equal to $g![\JJ]$, so
$[\JJ]\in\TCH^*(\JJ)_{\Q}$. Next, let us check that $\TCH^*(\JJ)_{\Q}$ is closed under the
Fourier transform $F:\CH^*(\JJ)_{\Q}\to\CH^*(\JJ)_{\Q}$. We have seen in the proof of Lemma
\ref{four-lem} that
$$F([n]_*\iota_*a)=\sum_{k\ge 0}\frac{n^k}{k!}\tau_k(a).$$
Since $F(x*y)=F(x)\cdot F(y)$, we derive the formula
\begin{equation}\label{F-prod-eq}
F(([n_1]_*\iota_*a_1)*\ldots*([n_s]_*\iota_*a_s))=
\sum_{k_1,\ldots,k_s}\frac{n_1^{k_1}\ldots n_s^{k_s}}{k_1!\ldots k_s!}\tau_{k_1}(a_1)\cdot\ldots\cdot
\tau_{k_s}(a_s).
\end{equation}
It remains to note that 
$$\tau_{k_1}(a_1)\cdot\ldots\tau_{k_s}(a_s)=T_{k_1}(0,a_1)\ldots T_{k_s}(0,a_s)([\JJ]),$$
hence, it is tautological by (ii). Thus, $\TCH^*(\JJ)_{\Q}$ is closed under the Fourier transform.
It follows that it is also closed under the usual product. Note that it contains all the classes
$\tau_k(a)=T_k(a)([\JJ])$. Now the fact that it is generated by these classes with respect to the usual product follows from \eqref{F-prod-eq}.

\noindent
(iv) Since the map $\si_*:\CH^*(\CC^{[\bullet]})\to\CH^*(\JJ)$ respects the Pontryagin products,
the first assertion follows from the formula
$$\si_{n*}\De_{n*}(a)=[n]_*\iota_*(a).$$
To check the assertion about the pull-backs we use the fact that the operators
$T_k(0,a)$ on $\CH^*(\CC^{[\bullet]})$ are multiplications by the pull-backs of $\tau_k(a)$
(see \eqref{T-k-0-eq}).
Since $[\CC^{[N]}]\in\TCH^0(\CC^{[N]})_{\Q}$ and the operators $T_k(0,a)$ preserve 
$\TCH^*(\CC^{[\bullet]})_{\Q}$, the result now follows from the fact that $\TCH^*(\JJ)_{\Q}$ is generated
by the classes $(\tau_k(a))$ with respect to the usual product (see part (iii)).

\noindent
(v) We can repeat the same proofs changing Chow groups to cohomology (and inserting 
appropriate signs where needed).
\ed

\begin{rems}
1. Consider the situation of Proposition \ref{diff-op-prop}(i). Assume that the algebra $A$ is equipped with a nonnegative $\Z$-grading (compatible with the $\Z/2\Z$-grading) such that  $\deg({\bf a}_0)=2$. Let us denote by $\DD'(A,{\bf a}_0)\sub\DD(A,{\bf a}_0)$ the subalgebra generated
by $\bP_{m,k}(a)$ with $m+k+\deg(a)\ge 2$. Then the analogue of Proposition \ref{diff-op-prop}(i)
holds if we only have an action of the subalgebra $\DD'(A,{\bf a}_0)$. In particular, this can be applied to
the actions of the Lie algebra $\HV'(\Z)\sub\HV(\Z)$ (see the Introduction). In the case of the action of
$\HV'(\Z)$ on $\CH^*(J)_{\Q}$ considered in \cite{P-lie} we recover the differential operators 
obtained in {\it loc. cit.}.

\noindent
2. The subalgebra $\DD_+(\CH^*(\CC),K)\sub\DD(\CH^*(\CC),K)$ 
considered in Proposition \ref{diff-op-prop}
is closely related to the algebra generated by the operators $(T_k(m,a))$ (see section \ref{Jac-sec}).
We will study these algebras in a sequel to this paper.
\end{rems}

Let us give an example of the calculation involving tautological cycles and using the operators
introduced in this paper.

\begin{prop}\label{pull-back-prop} 
(i) Consider the cycles $\tau_k(\CC)\in\CH^{k-1}(\JJ)$ (see \eqref{tau-eq}).
Their pull-backs under the morphisms $\si_N:\CC^{[N]}\to\JJ$ are given by
\begin{align*}
&\si_N^*\tau_k(\CC)=(-1)^kN^k\psi^{k-1}[\CC^{[N]}]+\\
&\sum_{i+n+m+p+l=k}(-1)^{n+l+m}\frac{k!}{(i+n)!(m+p)!l!}S(i+n,i)S(m+p,m)N^l\psi^{p+l}
\De_{i*}(K^n)*[p_0]^{*m}*[\CC^{[N-m-i]}]\\
&-\sum_{i+n+m+p+l+q=k, q\ge 1}(-1)^{n+l+m+q}\frac{k!q!}{(i+n+q)!(m+p)!l!}{i+q\choose i}
{m\choose q} S(i+n+q,i+q)S(m+p,m)N^l\times\\
&\psi^{p+l+n+q-1}[p_0]^{*(m+i)}*[\CC^{[N-m-i]}],
\end{align*}
where $[p_0]^{*m}$ denotes the $m$-th power with respect to the Pontryagin product,
$S(\cdot,\cdot)$ are the Stirling numbers of the second kind. Note that $\De_0:\CC\to S$
is just the projection to the base.

\noindent
(ii) For $k>g+\dim S+1$ one has $\tau_k(\CC)=0$, and hence, $\si_N^*\tau_k(\CC)=0$.
On the other hand, for $k>2g$ the class $\tau_k(\CC)$ becomes zero in $\CH^*(\JJ)_{\Q}$. 
\end{prop}

\Pf . (i) The idea is to use the formula
$$\si_N^*\tau_k(\CC)=T_k(0,\CC)([\CC^{[N]}])$$
(see \eqref{T-k-0-eq}). From Proposition \ref{T-P-prop} 
we get the following expression for $T_k(0,\CC)$:
$$T_k(0,\CC)=(-1)^kP_{1,1}(\CC)^k\psi^{k-1}+\sum_{i+n+j+l=k}(-1)^{n+j+l}\frac{k!}{(i+n)!j!l!}S(i+n,i)
P_{i,i}(K^n)P_{1,1}([p_0])^jP_{1,1}(\CC)^l\psi^l.$$
Recall that $P_{1,1}(\CC)x=Nx$ for $x\in\CH^*(\CC^{[N]})$.
Set $t=[p_0]\in\CH^1(\CC)$, $u=[\CC]\in\CH^0(\CC)$.
It follows easily from Proposition \ref{diff-op-prop}(ii) that the subalgebra
$\CH^*(S)[t,u^{[\bullet]}]\sub\wt{\TCH}^*(\CC^{[\bullet]})$ is preserved by
the operator $P_{1,1}([p_0])$, and its action is given by
$$P_{1,1}([p_0])|_{\Z[t,u^{[\bullet]}]}=t(\pa_u-\psi \pa_t).$$
From this we deduce by induction in $j$ that
$$P_{1,1}([p_0])^j([\CC^{[N]}])=P_{1,1}([p_0])^j(u^{[N]})=\sum_{m=0}^j(-\psi)^{j-m}S(j,m)t^mu^{[N-m]}$$
for $j\ge 0$.
Next, from Proposition \ref{diff-op-prop}(ii) we see that for an element $f\in\CH^*(S)[t,u^{[\bullet]}]$
one has
$$P_{i,i}(K^n)f=\De_{i_*}(K^n)\pa_u^i f-
t^i\cdot \sum_{q\ge 1}(-1)^q{i\choose q}\psi^{n+q-1}\pa^{i-q}_u\pa^q_t f.$$
From this we get formulas for $P_{i,i}(K^n)P_{1,1}([p_0])^j([\CC^{[N]})$ for all $i,j,n$. 
It remains to substitute them
into our expression for $T_k(0,\CC)([\CC^{[N]}])$.

\noindent (ii) The first assertion is clear since $\tau_k(\CC)\in\CH^{k-1}(\JJ)$. The second follows from the equality $\tau_k(\CC)=T_k(0,\CC)[\JJ]$ together with the fact that 
$T_k(0,\CC)=X_{0,k}$ is zero as an operator on $\CH^*(\JJ)_{\Q}$ for $k>2g$.
\ed

\begin{cor}\label{pull-back-cor} Assume that $S$ is a point. 

\noindent
(i) For $k>1$ one has 
\begin{align*}
&\si_N^*\tau_k(C)=\Ga_{e,k}*[C^{[N-k]}] 
-\sum_{i+m=k-1}(-1)^m{k\choose m}{i+1\choose 2}\De_{i*}(K)*[p_0]^{*m}*[C^{[N-k+1]}]\\
&-2\de_{k,2}[p_0]*[C^{[N-1]}],
\end{align*}
where
$$\Ga_{e,k}=\sum_{i+m=k}(-1)^m{k\choose m}[\De_i(C)]*[p_0]^{*m}\in\CH^{k-1}(C^{[k]})$$

\noindent
(ii) Assume that $K=(2g-2)[p_0]$. Then for $k>1$ one has
$$\si_N^*\tau_k(C)=\Ga_{e,k}*[C^{[N-k]}]-2g\de_{k,2}[p_0]*[C^{[N-1]}].$$
Without this assumption one has
$$\si_N^*\tau_k(C)\sim_{a.e.}\Ga_{e,k}*[C^{[N-k]}]-2g\de_{k,2}[p_0]*[C^{[N-1]}],$$
where $\sim_{a.e.}$ denotes algebraic equivalence.

\noindent
(iii) Assume that $K=(2g-2)[p_0]$. Then one has $\Ga_{e,k}=0$ for $k>g+1$. 
Also, $\Ga_{e,k}$ is a torsion class for $k\ge g/2+2$.

\noindent
(iv) One has $\Ga_{e,k}\sim_{a.e.} 0$ for $k>g+1$.
Furthermore, if $C$ admits a morphism of degree $d\ge 1$ to $\P^1$ 
then for $k\ge d+1$ some multiple of $\Ga_{e,k}$ is algebraically equivalent to zero.
\end{cor}

\Pf . (i) This follows from Proposition \ref{pull-back-prop}(i).

\noindent (ii) This follows from (i) since $\De_{i*}[p_0]=[p_0]^{*i}$.

\noindent (iii) Note that under these assumptions we have either $k>2$ or $g=0$. In either case
taking $N=k$ in (ii) we get $\si_k^*\tau_k(C)=\Ga_{e,k}$.
Now the first assertion follows from the trivial vanishing $\tau_k(C)=0$ for $k>g+1$. 
The second assertion follows from vanishing of $\tau_k(C)$ in $\CH^*(J)_{\Q}$ for $k\ge g/2+2$ that
is checked as follows. Note that $\tau_k(C)/k!=p_{k-1}$, where $(p_n)$ are the tautological
classes considered in \cite{P-lie}. Now the required vanishing follows from Proposition 4.2 of \cite{P-lie}
(note that the condition $K=(2g-2)[p_0]$ implies the vanishing of all the classes $q_n$).

\noindent (iv) The first statement follows from $K\sim_{a.e.}(2g-2)[p_0]$ as in (iii). For the second statement we use the fact that for $k\ge d+1$ some multiple of
$\tau_k(C)$ is algebraically equivalent to zero in this case
by the result of Colombo and Van Geemen in \cite{CG}.
\ed

\begin{rem}
The {\it modified diagonal classes} $\Ga_{e,k}\in\CH^{k-1}(C^{[k]})$ were introduced 
by Gross and Schoen in \cite{GS}. Parts (iii) and (iv) of the above
corollary are related to Propositions 4.5 and 4.8 of {\it loc. cit.}. In the case when $C$ is hyperelliptic 
and $p_0$ is stable under the hyperelliptic involution
Gross and Schoen show that $\Ga_{e,3}=0$. By
Corollary \ref{pull-back-cor}(ii) this is equivalent to the fact that $\tau_3(C)=0$ in this case 
(since $\si_N^*$ is injective for sufficiently large $N$). 
\end{rem}

\end{document}